\documentclass[11pt,a4paper]{amsart}
\usepackage{preamble}
\usepackage{ifsym}
\title{On Picard Groups of Perfectoid Covers of Toric Varieties}
\author[Gabriel Dorfsman-Hopkins]{Gabriel Dorfsman-Hopkins$^\dagger$}
\thanks{$^\dagger$UC Berkeley, 970 Evans Hall, Berkeley, CA 94720, United States.\\ \hspace*{0.6cm} \textit{gabrieldh@berkeley.edu}}

\author[Anwesh Ray]{Anwesh Ray$^\ddag$}
\thanks{$^\ddag$University of British Columbia, Vancouver, BC, V6T 1Z2, Canada.\\
\hspace*{0.6cm} \textit{anweshray@math.ubc.ca}}

\author[Peter Wear]{Peter Wear$^\star$}
\thanks{$^\star$University of Utah, 155 1400 E, Salt Lake City, UT 84112, United States.\\ \hspace*{0.6cm} \textit{wear@math.utah.com}}

\begin{document}


\begin{abstract} Let $X$ be a proper smooth toric variety over a perfectoid field of prime residue characteristic $p$.  We study the perfectoid space $\X$ which covers $X$ constructed by Scholze, showing that $\Pic(\X)$ is canonically isomorphic to $\Pic(X)[p^{-1}]$.  We also compute the cohomology of line bundles on $\X$ and establish analogs of Demazure and Batyrev-Borisov vanishing.  This generalizes the first author's analogous results for \textit{projectivoid space}.
\end{abstract}
\maketitle


\section{Introduction}

Perfectoid spaces are certain infinitely ramified nonarchimedean analytic spaces introduced by Scholze \cite{ScholzePS}. They have played a crucial role in settling a number of conjectures in arithmetic geometry and led to advances in $p$-adic Hodge theory, the Langlands program, and the study of Shimura varieties \cite{CaraianiScholze}, \cite{ScholzePadicHodgeTheory}, not least because they facilitate a correspondence between characteristic zero objects and their positive characteristic analogues. Given the paramount importance of perfectoid spaces it is of interest to develop some algebro-geometric tools to work with these highly non-noetherian objects.

In positive characteristic, perfectoid spaces arise functorially from varieties over $\F_p(t)$ by taking the \textit{completed perfection}: first taking the perfect closure and then completing the underlying rings with respect to the $t$-adic topology.  In mixed characteristic, the constructions are more subtle, but there are several known analogues including constructions for projective space, toric varieties \cite[Section 8]{ScholzePS}, and abelian varieties \cite{Peteretal}.

In \cite{DHthesis}, the first author studies the first case, considering the perfectoid analogue of projective space (dubbed \textit{projectivoid space}), and shows among other things that it has Picard group canonically isomorphic to $\Z[p^{-1}]$.  In this manuscript, we investigate whether the methods developed in \cite{DHthesis} apply to other classes of naturally arising perfectoid spaces, noticing that the construction of projectivoid space is a special case of the more general construction of perfectoid covers of toric varieties \cite[Section 8]{ScholzePS}. In positive characteristic, this construction coincides with the completed perfection. The following is a generalization of the first author's result \cite[Theorem 3.4]{DHthesis}.
\vspace{2cm}
\begin{Th}\label{main}
	Let $K$ be a perfectoid field with residue characteristic $p>0$, and $\Sigma$ be a complete smooth fan consisting of strongly convex rational cones. Let $X = X_{\Sigma,K}$ be the toric variety over $K$ associated to $\Sigma$ and let $\X\to X$ be its perfectoid cover.  Then the Picard group $\Pic(\X)$ is canonically isomorphic to $\Pic(X)[p^{-1}]$.
\end{Th}
We prove this in a similar way to \cite{DHthesis}, noticing first that $\X$ has a natural integral model whose special fiber is the scheme theoretic perfect closure of the toric variety associated to $\Sigma$ over the residue field of $K$.  The analogous result for scheme theoretic perfect closures is easily deduced, leaving two steps: first that every line bundle over the residue deforms uniquely to the integral model, and second that every line bundle on $\X$ extens uniquely to the integral model.

With Theorem \ref{main} we may identify a line bundle $\cL$ on $\X$ with a formal $p$-power root $\cM^{1/p^k}$ of a line bundle $\cM$ on $X$. We then compute the cohomology of $\cL$ in terms of the cohomology of powers of $\cM$, generalizing {\cite[Theorem 3.26]{DHthesis}}.

\begin{Th}\label{subMain}
	For every $i\ge0$, there are topologies for which the cohomology group $\HH^i(\X,\cL)$ is canonically isomorphic to the completion of $\colim\HH^i(X,\cM^{p^n})$.
\end{Th}

This allows us to deduce cohomological results on $\X$ from known results about the cohomology of toric varieties, and allows us to promote to the perfectoid setting the Demazure vanishing theorem \cite{Demazure} and the Batyrev-Borisov vanishing theorem \cite{BB}.

The authors point out that Heuer generalizes our results in \cite[Theorem 1.5]{Ben}, using the $v$-topology to study the Picard group of perfectoid covers of varieties with good reduction.

\subsection{Acknowledgments}
Dorfsman-Hopkins thanks Sebastian Bozlee and Martin Olsson for extremely enlightening comments and conversations.  Ray would like to thank Andres Fernandez Herrero, Kiran Kedlaya, Ravi Ramakrishna and Mike Stillman for some helpful conversations.  The authors also thank the anonymous referee for providing careful and detailed comments and suggestions, and for correcting an error in a previous draft. While working on this paper Dorfsman-Hopkins was partially supported by NSF grant DMS-1439786 while in residence at the Institute for Computational and Experimental Research in Mathematics in Providence, RI and NSF RTG grant DMS-1646385 as part of the Research Training Group in arithmetic geometry at the University of California, Berkeley.  Wear was partially supported by NSF grant DMS-1502651 and NSF RTG grant \#1840190.


\section{Toric Varieties and Toric Schemes}
We begin by reviewing the construction of a toric scheme associated to a fan and recalling a few useful properties.  One can find a complete reference about the theory over the complex numbers in \cite{CLS}, and over general fields in \cite{Danilov}.
\begin{Def}
	Let $N$ be a free abelian group of finite rank. Let $M := \Hom(N,\Z)$ be the dual lattice, and denote the canonical pairing $M_\R\times N_\R\to\R$ by $\langle m,n\rangle:= m(n)$.
	\begin{enumerate}[(i)]
		\item A \textit{strongly convex polyhedral rational cone} (henceforth, \textit{cone}) $\sigma \subseteq N_\R = N\otimes_\Z\R$ is a set of the form $\lambda_1\R_{\ge0} + \cdots + \lambda_t\R_{\ge_0}$ with $\lambda_i\in N$ subject to the condition that no line through the origin is contained in $\sigma$.  If the set $\lambda_1,\cdots,\lambda_t$ can be completed to a $\Z$-basis of $N$ then $\sigma$  is called \textit{smooth}.
		\item Given a cone $\sigma$ the \textit{dual cone} of $\sigma$ is
		\[\sigma^{\vee}:=\lbrace m\in M_{\R}\mid \langle m, u\rangle \geq 0\text{ for all }u\in \sigma\rbrace.\]
		\item A cone $\tau$ is said to be a \textit{face} of $\sigma$ if it is of the form $\sigma\cap H_m$ where to $m\in\sigma^\vee$, we associate the hyperplane $H_m=\{n\in N_\R : \la m,n\ra = 0\}$.
		\item A nonempty and finite collection of strongly convex polyhedral cones $\Sigma$ is called a \textit{fan} if it is closed under taking faces and if the intersection of any two cones is a face of both of them. A fan $\Sigma$ is \textit{smooth} if all of its cones are.
		\item The \textit{support} of $\Sigma$ is the union $|\Sigma|:=\bigcup_{\sigma\in\Sigma}\sigma$. If $|\Sigma|$ is all of $N_\R$, then $\Sigma$ is called \textit{complete}.
	\end{enumerate}
\end{Def}
The monoid $\sigma^{\vee}\cap M$ is finitely generated by \cite[Proposition 1.2.7]{CLS}.
\begin{Def}\label{toricDef}
	Let $A$ be a commutative ring and $\sigma\subseteq N_\R$ a strongly convex polyhedral cone. The affine toric scheme over $A$ associated to $\sigma$ is
	\[U_{\sigma,A} := \Spec A[\sigma^{\vee}\cap M].\]
	Each element $m\in\sigma^\vee\cap M$ induces an element of the ring of regular functions of $U_{\sigma,A}$.  We denote this function by $\chi^m\in A[\sigma^\vee\cap M]$ and call it the \textit{character} associated to $m$.
	When $A$ is understood we omit it from the notation.
\end{Def}
When $\sigma$ is smooth the structure of $U_{\sigma,A}$ is rather simple.
\begin{Ex}{\cite[Example 1.2.21]{CLS}}\label{structure}
	Let $A$ be a commutative ring and $\sigma$ a smooth strongly convex polyhedral cone.  Then there is a $\Z$-basis $e_1,\dots, e_n$ of $M$ such that $\sigma^{\vee}\cap M$ is generated by \[\{e_1,\cdots,e_r,\pm e_{r+1},\cdots,\pm e_n\}.\]In particular, there is an isomorphism
	\[A[\sigma^\vee\cap M]\cong A[x_1,\cdots,x_r,x_{r+1}^{\pm1},\cdots,x_n^{\pm1}]\]
	so that $U_{\sigma,A}$ is isomorphic to $\mathbb{A}^r_A\times\mathbb{G}_{m,A}^{n-r}$.
\end{Ex}
If $\tau$ is a face of $\sigma$ then the inclusion $\sigma^{\vee}\subseteq\tau^{\vee}$ induces a localization map $A[\sigma^\vee\cap M]\to A[\tau^\vee\cap M]$.  Passing to prime spectra induces an open immersion $U_\tau\hookrightarrow U_\sigma$.  Gluing along these open immersions allows us to define the toric scheme over $A$ associated to a fan $\Sigma$, which we denote by $X_{\Sigma,A}$.  As above, when $A$ is understood we omit it from the notation.  Notice also that this construction is compatible with base change in $A$.
\begin{Lemma}\label{toricBaseChange}
	If $\Sigma$ is a fan and $A\to B$ is a ring homomorphism, then there is a canonical isomorphism
	\[X_{\Sigma,A}\times_{\Spec A}\Spec B\cong X_{\Sigma,B}.\]
\end{Lemma}
\begin{proof}
	The affine case is clear, and since the isomorphisms are canonical over the affine cover, they glue.
\end{proof}
\begin{Remark}\label{relativeToricScheme}
	Given a fan $\Sigma$ and a base scheme $S$, one can define the \textit{relative toric scheme} $X_{\Sigma,S}\to S$ by building toric schemes as in Definition \ref{toricDef} affine locally on $S$, and then gluing canonically via Lemma \ref{toricBaseChange}.
\end{Remark}
When $A=k$ is a field, there is a well known structural result.
\begin{Def}
	A toric variety over $k$ is a normal separated scheme of finite type with an action by a split torus $T = \mathbb{G}_m^n$ such that there exists a point $x\in X(k)$ with trivial stabilizer and open dense orbit.  By \cite[Chapter 1, Theorem 6]{toroidal} every toric variety is isomorphic to $X_{\Sigma,k}$ for a unique fan $\Sigma$, and conversely if $\Sigma$ is a fan then $X_{\Sigma,k}$ is a toric variety.
\end{Def}
\begin{Remark}\label{alwaysNormal}
	In the literature one encounters toric varieties which may not be normal.  The toric varieties we study all correspond to fans, and are therefore normal by \cite[Chapter 1, Theorem 7]{toroidal}.  Therefore when we say toric variety we shall always assume normality.
\end{Remark}
We record a few more results that we will be using throughout the paper.
\begin{Prop}\label{proper}
	A relative toric scheme $X_{\Sigma,S}\to S$  is proper if and only if $\Sigma$ is complete.
\end{Prop}
\begin{proof}
	The case where $S=\Spec A$ is affine is \cite[Section 4 Proposition 4]{Demazure}, and properness is local on the target.
\end{proof}
\textcolor{black}{
It is shown in \cite[Corollary 7.4]{Danilov} that proper toric varieties have acyclic structure sheaves, and cohomology and base change theorems extend this property to the relative setting, implying in particular that a proper toric variety over any ring has an acyclic structure sheaf.
\begin{Prop}\label{localRingAcyclic}
	Let $S$ be a scheme and $\Sigma$ a complete fan, and $f:X_{\Sigma,S}\to S$ the associated relative toric scheme.  Then $R^if_*\cO_{X_{\Sigma,S}}=0$ for all $i>0$.  In particular, if $S=\Spec A$ is affine, $\cO_{X_{\Sigma,A}}$ is acyclic.
\end{Prop}
\begin{proof}
	By Lemma \ref{toricBaseChange}, $X$ is the pullback to $S$ of the proper $\bZ$ scheme $X_{\Sigma,\bZ}\to\Spec\bZ$.  By \cite[Corollary 7.4]{Danilov}, for each prime $\fp$ of $\bZ$, the module $\cO_{X_{\Sigma,\bZ}}\otimes k(\fp) = \cO_{X_{\Sigma,k(\fp)}}$ is acyclic.  In particular, for $i>0$ the function $\fp\mapsto\dim_{k(\fp)}\HH^i(X_{k(\fp)},\cO_{X_{\Sigma,k(\fp)}})$ is the constant zero function on $\Spec\bZ$. Therefore by Grauert's theorem of cohomology and base change \cite[Theorem 28.1.5]{Vakil} $\cO_{X_{\Sigma,\bZ}}$ is acyclic and its cohomology commutes with arbitrary base change, whence the result follows.
\end{proof}
A toric variety $X_\Sigma$ comes equipped with an affine cover $\fU = \{U_{\sigma}\to X_\Sigma\}$ as $\sigma$ varies over the maximal cones of $\Sigma$.  By \cite[Theorem 2.1]{Gub}, every line bundle on $X_\Sigma$ trivializes on $\fU$ so that much of this study reduces to the combinatorics of $\Sigma$.
}
\begin{Prop}\label{toricPic}
	Let $\Sigma$ be a fan, $k$ a field, and $X = X_{\Sigma,k}$ the associated toric variety.  The Picard group of $X$ does not depend on the field $k$.
\end{Prop}
\begin{proof}
	Let $\{\sigma_1,\sigma_2,\cdots,\sigma_r\}\subset\Sigma$ be the set of maximal cones. Set $U_i = U_{\sigma_i,k}$, and denote by $\mathfrak{U}$ the cover $\{U_i\to X\}$.  The finite intersections of the $U_i$ are affine toric varieties and therefore by \cite[Theorem 2.1]{Gub} all have trivial Picard group.  Therefore the \v{C}ech-to-derived functor spectral sequence degenerates to an isomorphism
	\[\Pic(X) \cong \HH^1(X,\cO_X^*)\cong\cHH^1(\fU,\cO_X^*).\]
	Consider the first three terms of the \v{C}ech sequence,

 \vspace{0.25cm}\[C^0(\fU,\cO_X^*)\longtoo{\delta^0}C^1(\fU,\cO_X^*)\longtoo{\delta^1}C^2(\fU,\cO_X^*)\longto\cdots\]
 \vspace{0.25cm}
 \\
	Notice that for any cone $\sigma$, $\HH^0(U_\sigma,\cO_X)\subseteq k[x_1^{\pm1},\cdots,x_n^{\pm1}]$. Any $f\in \HH^0(U_\sigma,\cO_X^*)$ must therefore be a monomial $f = \lambda x_1^{m_1}\cdots x_n^{m_n} = \lambda\chi^m$ for some $m\in\sigma^\vee\cap M$ and $\lambda\in k^*$.  Let $D^1\subseteq\ker\delta^1$ be the subgroup of all cocycles consisting of monic monomials $(\chi^{m_{ij}})$, and let $D^0 = D^1\cap\im\delta^0$.  We have the following map of short exact sequences
  \vspace{0.5cm}

	\[
		\begin{tikzcd}
			0\rar & \im\delta^0\rar & \ker\delta^1\rar & \cHH^1(\fU,\cO_X^*)\rar& 0\\
			0\rar & D^0\rar\ar{u} & D^1\rar\ar{u} & \tilde{H}\ar{u}\rar& 0
		\end{tikzcd}
	\]
	\vspace{0.5cm}
  \\
  where $\tilde H$ is defined as the cokernel of the inclusion $D^0\into D^1$. As $D^0$ and $D^1$ consist only of monic monomials, $\tilde H$ does not depend on $k$, only the combinatorics of $\Sigma$.  Therefore, if we prove the vertical map on the right is an isomorphism, we will be done.

	Fix a cocycle $\alpha = (\alpha_{ij}) = (\lambda_{{ij}}\chi^{m_{ij}})\in\ker\delta^1$.  Notice that the class of $\alpha$ does not depend on the $\lambda_{ij}$.  Indeed, the class $(\lambda_{ij}) = \delta^0(\lambda_{1r},\lambda_{2r},\cdots,\lambda_{r-1,r},1)$ (noticing that the cocycle condition implies that $\lambda_{ir}/\lambda_{jr} = \lambda_{ij}$).  Therefore we may assume that $\alpha$ is in $D^1$ without changing the associated cohomology class, and so $\tilde H\to\cHH^1(\fU,\cO_X^*)$ surjects (by the commutativity of the diagram above).  Let $K$ be the kernel, and let $Q^0$ and $Q^1$ be the cokernels of the injections $D^0\to\im\delta^0$ and $D^1\to\ker\delta^1$ respectively.  Then the snake lemma gives us the following exact sequence
	\[0\longto K\longto Q^0\longto Q^1\longto 0.\]
	If we show the map on the right injects we are done.  But
	\[Q^0 = \{(\lambda_i\lambda_j^{-1}\chi^{m_i-m_j})\}/\{(\chi^{m_i-m_j})\} = \{(\lambda_i/\lambda_j)\} \cong (k^*)^{(r-1)}.\]
	Furthermore, since the $\lambda_{ij}\chi^{m_{ij}}$ must satisfy the cocycle condition, we can deduce the $\lambda_{ij}$ once we know them for one fixed $j$, so that
	 \[Q^1 = \{\lambda_{ij}\chi^{m_{ij}}\}/\{\chi^{m_{ij}}\} \cong (k^*)^{(r-1)}.\]
	This identifies $Q^0\to Q^1$ with the identity map and so we are done.
\end{proof}
The proof of Proposition \ref{toricPic} has a useful corollary which we record here for later use.  We will first need a definition.
\begin{Def}[The $n$th power map]\label{powerMap}
	Fix a fan $\Sigma$ and an integer $n$.  For each $\sigma\in\Sigma$ the multiplication by $n$ on $M$ map induces a homomorphism $A[M\cap\sigma^\vee]\to A[M\cap\sigma^\vee]$.  This map is compatible with restrictions to faces $\tau\le\sigma$, and therefore glues to a morphism $\varphi_n:X_{\Sigma,A}\to X_{\Sigma,A}$, which we will call the \textit{$n$th power map on $X_{\Sigma,A}$}
\end{Def}
\begin{Cor}\label{powerPullback}
	Let $k$ be a field and $\Sigma$ a fan.  The $n$th power map on $X =X_{\Sigma,k}$ induces a pullback homomorphism $\varphi_{n}^*:\Pic X\to\Pic X$ which can be identified with multiplication by $n$.
\end{Cor}
\begin{proof}
	Let $\fU$ be the same cover from the proof of Proposition \ref{toricPic}, and recall that an arbitrary \v{C}ech cocycle in $\cHH^1(\fU,\cO_X^*)$ may be represented by monic monomials, thus be of the form $(\chi^{m_{ij}})$.  Then pullback along $\varphi_n$ corresponds to multiplication by $n$ on $M$ so that it takes $(\chi^{m_{ij}})\mapsto(\chi^{nm_{ij}})$.  But this is precisely the multiplication by $n$ map on (monic monomial) cocycles, so that it descends to the same on cohomology.
\end{proof}

\section{Perfectoid Covers of Toric Varieties}\label{PCoverSection}
To begin this section we review Scholze's construction of the perfectoid cover of a toric variety $X=X_{\Sigma,K}$. We use Huber's theory of \textit{adic spaces}, a non-archimedean analogue of complex analytic spaces, which are general enough to handle the infinite ramification required to define perfectoid spaces.  We will not develop the entire theory of adic spaces here, but refer the reader to \cite{Huber1} and \cite{Huber2} where it was introduced, or to \cite[Section 2]{ScholzePS} for an introduction with more of an emphasis on perfectoid spaces. We do introduce a few preliminary notions.
\begin{Def}\label{tateDef}
	A \textit{Tate ring} is a complete topological ring $A$ such that there exists a topologically nilpotent unit called the \textit{pseudouniformizer} $\varpi\in A$ and an open subring $A_0\subset A$ containing $\varpi$ such that $A_0$ has the $\varpi$-adic topology. A ring of integral elements $A^+\subseteq A$ is an open integrally closed subring of the power bounded elements $A^\circ$.
\end{Def}
Associated to such a pair $(A,A^+)$ is the \textit{adic spectrum} $\Spa(A,A^+)$ which is a topological space whose underlying set consists of equivalence classes of valuations on $A$ which are bounded by $1$ on $A^+$.  If $x\in\Spa(A,A^+)$ is a point, we suggestively denote the associated valuation by $f\mapsto|f(x)|$ for $f\in A$.  General adic spaces are locally of the form $\Spa(A,A^+)$, and all adic spaces come equipped with various sheaves related to their structure sheaf.  We make heavy use of these sheaves so we define them carefully.
\begin{Def}
Let $\cX$ be an adic space with structure sheaf $\cO_\cX$.  \textit{The sheaf of integral elements} $\cO_\cX^+$, a sheaf of (topological) rings defined by the rule
\vspace{0.25cm}
\[\cO_\cX^+(U)=\left\{f\in\cO_\cX(U)\text{ }:\text{ }|f(x)|\le 1\text{ for all }x\in U\right\}.\]
\vspace{0.25cm}
\\
The sheaf of integral elements has an ideal $\cO_\cX^{++}\subseteq\cO_\cX^+$ called \textit{the sheaf of topologically nilpotent elements} defined by the rule
\vspace{0.25cm}
\[\cO_\cX^{++}(U)=\left\{f\in\cO^+_\cX(U)\text{ }:\text{ }|f(x)|< 1\text{ for all }x\in U\right\}.\]
\vspace{0.25cm}
\\

\end{Def}
Perfectoid spaces are adic spaces built from special kinds of rings.
\begin{Def}\label{perfectoidDef}
	A topological ring $A$ is a \textit{perfectoid Tate ring} if it is a Tate ring containing a pseudouniformizer $\varpi\in A$ such that
	\begin{itemize}
		\item $A$ is \textit{uniform}, i.e. $A^\circ$ is bounded,
		\item $\varpi^p|p$ in $A^\circ$,
		\item the $p$-power Frobenius map is an isomorphism
		      \[\Phi:A^\circ/\varpi \longtoo{\sim} A^\circ/\varpi^p.\]
	\end{itemize}
	A (Tate) perfectoid space is an adic space which is locally the adic spectrum of perfectoid Tate rings.
\end{Def}
\begin{Remark}
	There is a more general framework of so called \textit{Huber rings} and \textit{perfectoid rings} from which one can construct adic spaces and perfectoid spaces which are not Tate.  The main difference here is that one does not assert the existence of a pseudouniformizer as in Definitions \ref{tateDef} and \ref{perfectoidDef}.  This allows one to apply the theory to rings which are more integral in nature.  Nevertheless, we will not require this level of generality because we will always be working over a perfectoid base field and will therefore inherit a pseudouniformizer from this field.
\end{Remark}
Let $K$ be a perfectoid field with valuation ring $K^\circ$ and residue characteristic $p$.  To any cone $\sigma$ we can associate the Tate ring $K\langle\sigma^\vee\cap M\rangle$ of convergent power series in the characters of $\sigma$.  We form the associated adic space
\[\mathcal U_\sigma: = \Spa\left(K\langle\sigma^\vee\cap M\rangle,K^\circ\langle\sigma^\vee\cap M\rangle\right).\]
As with toric spaces, if $\tau$ is a face of $\sigma$ we have an induced open immersion $\cU_\tau\subset\cU_\sigma$ (we note that this is generally open in the analytic topology rather than the much coarser Zariski topology).  Therefore, if $\Sigma$ is a fan we can form an adic space $\cX_\Sigma$.
\begin{Remark}\label{adicAnalytification}
	There is a fully faithful analytification functor $X\mapsto X^{ad}$ from schemes over $K$ to adic spaces over $K$, but if $X_{\Sigma}$ is not proper then $\cX_\Sigma$ is not in general isomorphic to $(X_\Sigma)^{ad}$.  For example, if $\Sigma$ is the cone associated to $\mathbb{N}^n\subseteq\Z^n$, then $\cX_\Sigma$ is the rigid analytic disk whereas $(X_\Sigma)^{ad}$ is affine space.  Nevertheless, if $\Sigma$ is a complete fan (so that $X_\Sigma$ is proper), then we do have $(\cX_\Sigma)=(X_\Sigma)^{ad}$ by \cite[Section 8]{ScholzePS}.
\end{Remark}
As in Definition \ref{powerMap}, multiplication by $p$ on $M$ induces a map $K\langle\sigma^\vee\cap M\rangle\to K\langle\sigma^\vee\cap M\rangle$ inducing the $p$-power map $\varphi$ on the adic spectra.  Passing to the completed direct limit results in the perfectoid algebra $K\langle\sigma^\vee\cap M[p^{-1}]\rangle$ and on adic spectra this translates to a so called \textit{tilde inverse limit} (as in \cite{ScholzePS}):
\[\cU_\sigma^\perf\sim\ilim_{\varphi}\cU_\sigma.\]
This construction globalizes.  If $\tau$ is a face of $\sigma$, then $K\la\sigma^\vee\cap M[p^{-1}]\ra\to K\la\tau^\vee\cap M[p^{-1}]\ra$ is a rational localization so that $\cU_\tau^{\perf}\to \cU_\sigma^{\perf}$ is an open immersion.  Therefore given a fan $\Sigma$, we can glue the affinoid perfectoid spaces $\cU_\sigma^\perf$ along their intersections and construct a perfectoid space $\cX_\Sigma^\perf$.  As in Definition \ref{powerMap}, the map $\varphi$ commutes with this gluing, so that it induces a map of adic spaces $\cX^\perf_\Sigma\to\cX_\Sigma$, and checking locally we see that in fact $\cX^\perf_\Sigma\sim\ilim_\varphi\cX_\Sigma$.
\begin{Def}\label{perfectoidCover}
Given a complete fan $\Sigma$ the map $\cX^\perf_\Sigma\to\cX_\Sigma$ is called the \textit{perfectoid cover} of the toric variety $\cX_\Sigma$.
\end{Def}
If $K$ has characteristic $p>0$ then $\varphi$ can be identified with relative Frobenius over $K$ and $\cX_\Sigma^\perf$ is the completed perfection of the adic space $\cX_\Sigma$.
\begin{Remark}
	Inverse limits do not exist in general in the category of adic spaces.  That being said, if affinoid locally the completed direct limit of the ring of functions is perfectoid, then the perfectoid space built from those rings satisfies the universal property of inverse limit among all perfectoid spaces.  It is in this sense that the tilde inverse limit is well defined.

	The tilde inverse limit agrees with a categorical inverse limit if we enlarge our category.  More precisely, perfectoid spaces embed in the category of \textit{diamonds} \cite{ScholzeDiamonds} which do have inverse limits.  In this case the diamond in the limit is representable by a perfectoid space, and this perfectoid space is the tilde inverse limit.  We will not use this formalism.
\end{Remark}
\begin{Ex}\label{examples}
	Let us see how this works in a couple of examples.
	\begin{enumerate}
		\item Let $\sigma$ be the cone associated to $\mathbb{N}^n\subseteq\mathbb{Z}^n$.  Then $K\la\sigma^\vee\cap M\ra = K\la x_1,\cdots,x_n\ra$ is the Tate algebra, and the associated adic space $\cU_\sigma$ is the rigid analytic unit disk $\mathbb{D}^{n}$.  The map
		\[\varphi:K\la x_1,\cdots,x_n\ra\to K\la x_1,\cdots,x_n\ra\]
		is the homomorphism taking $x_i\mapsto x_i^p$ for each $i$, and the completed direct limit along $\varphi$ is the perfectoid Tate algebra
		\[K\la\sigma^\vee\cap M[p^{-1}]\ra = K\la x_1^{1/p^\infty},\cdots,x_n^{1/p^\infty}\ra.\]
		Thus the perfectoid cover $\cU_\sigma^\perf$ is the perfectoid unit disk $\mathbb{D}_n^\perf$.
		\item Let $\sigma$ be the cone associated to $\{0\}$. Then $K\la\sigma^\vee\cap M\ra = K\la x_1^{\pm1},\cdots,x_n^{\pm1}\ra$ and $\cU_\sigma$ is the rigid analytic torus $\bT^n$. $\cU_\sigma^\perf$ is the perfectoid torus $\bT^{n,\perf}$ associated to
		\[K\left\la\sigma^\vee\cap M[p^{-1}]\right\ra = K\la x_1^{\pm1/p^\infty},\cdots,x_n^{\pm1/p^\infty}\ra.\]
		\item Let $\sigma$ be a smooth cone.  As in Example \ref{structure} we may take $\sigma^\vee\cap M$ to be generated by $\{e_1,e_2,\cdots,e_r,\pm e_{r+1},\cdots,\pm e_n\}$.  Then
		\[K\la\sigma^\vee\cap M\ra = K\la x_1,\cdots,x_r,x_{r+1}^{\pm1},\cdots,x_n^{\pm1}\ra,\]
		so that $\cU_\sigma \cong \mathbb D^{r}\times\bT^{n-r}$.  Taking limits along $\varphi$ we have
		\[K\la\sigma^\vee\cap M[p^{-1}]\ra \cong K\la x_1^{1/p^\infty},\cdots,x_r^{1/p^\infty},x_{r+1}^{\pm1/p^\infty},\cdots,x_n^{\pm1/p^\infty}\ra,\]
		so that $\cU_\sigma^{\perf} \cong \mathbb D^{r,\perf}\times\bT^{n-r,\perf}$.
		\item \textcolor{black}{Let $e_1,\cdots,e_n$ be a basis for the free abelian group $N$, and let $e_0 = -e_1-e_2-\cdots-e_n$.  Define $\Sigma$ to be the fan generated by proper subsets of $\{e_0,\cdots,e_n\}$.  Then one can check \cite[Example 3.1.10]{CLS} that  $X_\Sigma = \bP^n$.  As this is proper, the discussion in Remark \ref{adicAnalytification} implies that} the adic space associated to $\Sigma$ is also projective space.  One can identify the map $\varphi$ with the map
    \[[x_0:\cdots:x_n]\mapsto[x_0^p:\cdots:x_n^p]\]
    in projective coordinates, so that the associated perfectoid cover is \textit{projectivoid space} $\mathbb{P}^{n,\perf}$, which is studied extensively in \cite{DHthesis}.
	\end{enumerate}
\end{Ex}

Many of our computations will be working explicitly with various \v{C}ech complexes associated to the cover $\fU$ of $\cX^\perf$. They will all be built by taking units and quotients of the following complex.

\begin{Def}\label{intCech}
Let $\Sigma$ be a complete fan, and consider the associated perfectoid space $\X := \cX_\Sigma^\perf$.  Let $\{\sigma_1,\cdots,\sigma_r\}$ be the maximal cones of $\Sigma$ and consider the cover $\fU = \{\cU_{\sigma_i}^\perf\to\X\}$. We define $C_{\X}^\bullet:=\check C(\fU, \cO^+_{\cX^\perf})$ to be the \v{C}ech complex for the sheaf $\cO^+_{\cX^\perf}$ with respect to the cover $\fU$.  Concretely, we have \vspace{0.25cm}
  \begin{eqnarray*}
    C_{\X}^\bullet = 0\to\prod_i K^\circ\left\langle\sigma_i^\vee\cap M[p^{-1}]\right\rangle\to\prod_{i<j}K^\circ\left\langle(\sigma_i\cap\sigma_j)^\vee\cap M[p^{-1}]\right\rangle\to\cdots\\
    \cdots\to K^\circ\left\langle(\sigma_1\cap\cdots\cap\sigma_r)^\vee\cap M[p^{-1}]\right\rangle\to0.
  \end{eqnarray*}
  \vspace{0.25cm}
\end{Def}

\begin{Prop}\label{intAcyclic}
	 With notation as in Definition \ref{intCech}, we have $\HH^i(C_{\X}^\bullet)= 0.$
\end{Prop}
\begin{proof}
  We first remark that by construction we have $\cO_{\X}^+(\cU^\perf_\sigma) = K^\circ\la\sigma^\vee\cap M[p^{-1}]\ra$ for every $\sigma\in\Sigma$.  Note also that multiplication by $p$ induces an isomorphism on $M[p^{-1}]$, and therefore an isomorphism on $K^\circ\la\sigma^\vee\cap M[p^{-1}]\ra$.  This commutes with differentials on the \v{C}ech complex $C^\bullet_{\X}$, so that passing to cohomology we get an isomorphism on $\HH^i(C^\bullet_{\X})$.  We would like to use this to clear denominators on the exponents of our cocycles and reduce to the case of Proposition \ref{localRingAcyclic}, but in order to show a finite power of $p$ will clear these denominators, we have to show that it suffices to consider cocycles represented by polynomials.  To do so we construct the following cochain complex:
  \vspace{0.25cm}
  \begin{eqnarray*}
    C^\bullet_{\X, f} = 0\to\prod_i K^\circ\left[\sigma_i^\vee\cap M[p^{-1}]\right]\to\prod_{i<j}K^\circ\left[(\sigma_i\cap\sigma_j)^\vee\cap M[p^{-1}]\right]\to\cdots\\
    \cdots\to K^\circ\left[(\sigma_1\cap\cdots\cap\sigma_r)^\vee\cap M[p^{-1}]\right]\to0.
  \end{eqnarray*}
  \vspace{0.25cm}
  \\
Notice that the (termwise) $\varpi$-adic completion of this complex is $C^\bullet_{\X}$.  The terms in the sequence are linearly topologized $K^\circ$-modules with a countable system of fundamental neighborhoods of 0 (given by the $\varpi^n$), and the differentials are continuous maps.  Therefore, applying \cite[Tag 0AS0]{Stacks} in descending induction, the associated short exact sequences giving the cohomology of $C^\bullet_{\X,f}$ remain exact upon completion.  In particular, the completion of this sequence commutes with taking cohomology, so that we have that $H^i(C^\bullet_{\X}) = \HH^i(C^\bullet_{\X,f})^\wedge{}$.  Therefore it suffices to prove that the nonzero cohomology of $C^\bullet_{\X,f}$ vanishes.

Let $X_0 := X_{\Sigma,K^\circ}$ be the toric scheme associated to $\Sigma$ over $K^\circ$, and consider the induced affine cover $\fU_0 = \{U_{\sigma_i,K^\circ}\to X_0\}$.  As $X_0$ is proper, \v{C}ech cohomology with respect to an affine cover computes quasicoherent cohomology, so that applying Proposition \ref{localRingAcyclic} we have for $i>0$,
\[\cHH^i(\fU_0,\cO_{X_0}) = \HH^i(X_0,\cO_{X_0}) = 0.\]
Furthermore, the natural inclusions $K^\circ[\sigma^\vee\cap M]\into K^\circ[\sigma^\vee\cap M[p^{-1}]]$ for each $\sigma\in\Sigma$ induce the following map of cochain complexes:
\[
\begin{tikzcd}
  \cdots\rar & C^{i-1}(\fU_0,\cO_{X_0})\ar[r,"\delta_0^{i-1}"]\ar[d,hookrightarrow] & C^{i}(\fU_0,\cO_{X_0})\ar[r,"\delta_0^{i}"]\ar[d,hookrightarrow] & C^{i+1}(\fU_0,\cO_{X_0})\rar\ar[d,hookrightarrow] & \cdots\\
  \cdots\rar & C^{i-1}_{\X,f}\ar[r,"\delta^{i-1}"] & C^i_{\X,f}\ar[r,"\delta^{i}"] & C^{i+1}_{\X,f}\rar & \cdots
\end{tikzcd}
\]
As above, multiplication by $p$ on $M[p^{-1}]$ induces an isomorphism of $C^\bullet_{\X,f}$ which commutes with differentials.  In particular it induces an isomorphism of $\HH^i(C^\bullet_{\X,f})$ with itself.  Fix some $\alpha = (\sum \lambda_m \chi^m)_I\in C^i_{\X,f}$ where $I$ indexes over the product.  We denote by $p*\alpha = (\sum \lambda_m \chi^{pm})_I$.  If $\alpha\in\ker\delta^i$, then so is $p*\alpha$.  For some $N$ we have $p^N*\alpha\in C^i(\fU_0,\cO_{X_0})$, and thus in $\ker\delta^i_0$.  But $\cO_{X_0}$ is \v{C}ech acyclic so $p^N*\alpha=\delta_0^{i-1}(\beta)$ for some $\beta$.  Viewing $\beta$ as an element of $C^{i-1}_{\X,f}$ we see that
\[\delta^{i-1}(p^{-N}*\beta) = p^{-N}*\delta^{i-1}(\beta) = p^{-N}*p^N*\alpha = \alpha\]
In particular, $C^\bullet_{\X,f}$ is exact in the $i$th position and we win.

\end{proof}

\subsection{Triviality for Picard Groups for Affinoid Perfectoid Toric Spaces}\label{acyclicitySection}
In order to study the Picard group of $\X:=\X_\Sigma$, we would like to make use of the cover $\fU = \{\cU^\perf_\sigma\to\X_\Sigma\}$ induced by the maximal cones of $\sigma\in\Sigma$.  In particular, if every line bundle on $\X_\Sigma$ trivialized on this cover and its various intersections, it would suffice to study the \v{C}ech cohomology group $\cHH^1(\fU,\cO_{\X}^*)$, which is accessible in practice.  By a limiting argument it will suffice to show that the affinoid neighborhoods $\cU_\sigma$ of the toric \textit{variety} $\cX_\Sigma$ have trivial Picard groups.
\begin{Def}\label{fanInducedTrivialization}
	Let $\cX_\Sigma$ (respectively $\X_\Sigma$) be the rigid space (resp. perfectoid space) associated to a fan $\Sigma$.  If the Picard groups of the affinoid neighborhoods $\cU_\sigma$ (resp. $\cU_\sigma^\perf$) are trivial for each cone $\sigma\in\Sigma$, we say that the Picard group of $\cX_\Sigma$ (resp. $\X_\Sigma$) \textit{trivializes on the analytic affinoid cover induced by the fan}.
\end{Def}
\begin{Lemma}\label{smoothPic}
	If $\cX_\Sigma$ is a rigid space associated to a smooth fan, then the Picard group of $\cX_\Sigma$ trivializes on the analytic affinoid cover induced by the fan.
\end{Lemma}
\begin{proof}
	If $\sigma$ is a smooth cone then $\cU_\sigma = \bD^r\times\bT^{n-r}$ as in Example \ref{examples}(3).  Then the Picard group of $\cU_\sigma$ corresponds to projective modules on $K\la x_1,\cdots,x_r,x_{r+1}^{\pm1},\cdots,x_n^{\pm1}\ra$, which were shown to be trivial in \cite[Satz I]{Lutke}.
\end{proof}
We spend the rest of the section proving the following result.
\begin{Prop}\label{trivialPic}
	Let $\cX_\Sigma$ be a toric rigid space with perfectoid cover $\X_\Sigma$.  If the Picard group of $\cX_\Sigma$ trivializes on the analytic affinoid cover induced by the fan, then the same can be said for $\X_\Sigma$.
\end{Prop}
We will use a limiting argument together with the following theorem of Gabber and Ramero.
\begin{Prop}[{\cite[5.4.42]{GaRo}}]\label{gabber}
  Let R be a commutative ring, $t\in R$ a nonzero divisor, and $I\subseteq R$ an ideal. Let $\hat{R}$ be the $(t,I)$-adic completion of $R$, and suppose $(R, tI)$ form a henselian pair. Then the base extension functor $R[t^{-1}]\textbf{-Mod}\to\hat R[t^{-1}]\textbf{-Mod}$ induces a bijection between isomorphism classes of finite projective $R[t^{-1}]$-modules and finite projective $\hat R[t^{-1}]$-modules.
\end{Prop}
Let us formulate the specific consequence of Proposition \ref{gabber} in the way we will use it.
\begin{Cor}\label{gabberCorollary}
	Let $R = \cup R_i$ be the union of $t$-adically complete rings, and let $\hat R$ be its $t$-adic completion.  Then the natural map $\Pic R[t^{-1}]\to\Pic\hat R[t^{-1}]$ is an isomorphism.  Furthermore, if $\Pic R_i[t^{-1}]=0$ for all $i$, then $\Pic\hat R[t^{-1}]=0$
\end{Cor}
\begin{proof}
	We first establish the isomorphism $\Pic R[t^{-1}]\to\Pic\hat R[t^{-1}]$, for which it suffices to show that $(R,tR)$ form a henselian pair so that the result follows immediately from Proposition \ref{gabber}.  Suppose $f(x)\in R[x]$ is monic, and that after reducing mod $t$,  $\overline f(x) = g_0(x)h_0(x)$ with $g_0,h_0$ monic.  For some large $i$, we have $f(x)\in R_i[x]$.  Furthermore, one checks (using for example \cite[Tag 00DD]{Stacks}) that $R/tR\cong\colim R_i/tR_i$.  Therefore (perhaps increasing $i$) we also have $g_0(x),h_0(x)\in (R_i/tR_i)[x]$.  But as $R_i$ is $t$-complete, $(R_i,tR_i)$ form a henselian pair \cite[Tag 0ALJ]{Stacks}, and so the factorization lifts to $f(x)=g(x)h(x)\in R_i[x]\subseteq R[x]$.  Therefore $(R,tR)$ form a henselian pair, and we have established the desired isomorphism

	For the second statement it now suffices to show that every invertible $R[t^{-1}]$-module is free.  By \cite[Tag 0B8W]{Stacks}, every invertible $R[t^{-1}]$-module is the base extension of some invertible $R_i[t^{-1}]$-module, which is free by assumption.
\end{proof}
\begin{proof}[Proof of Proposition \ref{trivialPic}]
	Fix $\sigma\in\Sigma$ and let $R_i = K^\circ\la \sigma^\vee\cap p^{-i}M\ra$, and let $R = \cup R_i$.  Then $K\la\sigma^\vee\cap M[p^{-1}]\ra = \hat R[1/\varpi]$ where $\varpi$ is a pseudouniformizer of $K$.  Therefore by Corollary \ref{gabberCorollary} it suffices to show that every invertible $R_i[1/\varpi]$ module is free.  But as these are $K\la \sigma^\vee\cap p^{-i}M\ra\cong K\la\sigma^\vee\cap M\ra$, which have trivial Picard group is trivial by assumption.
\end{proof}
Lemma \ref{smoothPic} gives the following immediate corollary.
\begin{Cor}
	Let $\Sigma$ be a smooth cone.  Then the Picard group of $\X_\Sigma$ trivializes on the affinoid cover induced by the fan.
\end{Cor}
\begin{Remark}
	If $X_\Sigma$ is a toric variety (viewed as a scheme), the analogous statement for the Zariski affine cover of $X_\Sigma$ always holds by \cite[Theorem 2.1]{Gub}.  Therefore, if the Picard group of $K[\sigma^{\vee}\cap M]$ remains trivial after completion, all of the rigid toric spaces we are studying satisfy Definition \ref{fanInducedTrivialization}, and the smoothness hypothesis can be dropped.

	One might attack this in the following way: Let $k = K^\circ/K^{\circ\circ}$ be the special fiber.   Then \cite[Theorem 2.1]{Gub} implies that the Picard group of $k[\sigma^\vee\cap M]$ is trivial, and a simple deformation argument \cite[Corollary 2.15]{DHthesis} shows the same can be said for $K^\circ\la\sigma^\vee\cap M\ra$.  One then must argue that the Picard group remains trivial after inverting $\varpi$.
\end{Remark}

\section{Comparison Isomorphisms}
We now prove the main theorem of this manuscript.  For the rest {of this section} we fix a perfectoid field $K$ with valuation subring $K^\circ$, maximal ideal $\mathfrak{m}$, and residue field $k=K^\circ/\mathfrak{m}$.  Set $\Sigma$ to be a complete fan, and assume that the Picard group of the rigid space $\cX_\Sigma$ is trivialized on the affinoid cover induced by $\Sigma$ (for example if $\Sigma$ is smooth).  We fix two schemes $X = X_{\Sigma,K}$ and $X_0 = X_{\Sigma,k}$, and let ${\X} = \cX_{\Sigma,K}^\perf$ be the perfectoid cover of $X$. Also let $\{\sigma_1,\cdots,\sigma_r\}$ be the maximal cones of $\Sigma$, and let $\fU = \{\cU_{\sigma_i,K}^\perf\to{\X}\}$ be the induced cover of ${\X}$.  For simplicity, we denote by $\cU^{\perf}_i := \cU_{\sigma_i,K}^\perf$ and $\cU^\perf_{i_1\cdots i_t} := \cU^{\perf}_{\sigma_{i_1}\cap\cdots\cap\sigma_{i_t}} = \cU^{\perf}_{i_1}\cap\cdots\cap\cU^{\perf}_{i_t}$. Our main theorem is the following.
\begin{Th}\label{main2}
  There is a canonical isomorphism $\Pic({\X})\cong\Pic(X)[p^{-1}]$.
\end{Th}

{Our rough strategy consists of 3 main steps.
\begin{itemize}
	\item{\textbf{Extension Step:} In Section \ref{integral}, we show that each line bundle on $\X$ extends uniquely to the integral model. The inclusion $\cO_{\X}^{+*}\into\cO_{\X}^*$ induces a map $\HH^1(\X,\cO_{\X}^{+*})\to\HH^1(\X,\cO_{\X}^*)$, which we show is an isomorphism.}
	\item{\textbf{Deformation Step:} In Section \ref{deformation}, we show that each line bundle on the special fiber which is trivial on the cover $\fU$ deforms uniquely to an integral model. To do this, we take the \v{C}ech complex $C^{*\bullet}_{\X}$ of $\cO_{\X}^{+*}$ along the standard cover and compare it with its reduction modulo topological nilpotents $\overline{C^{*\bullet}_{\X}}$}.
	\item{\textbf{Explicit Comparison:} In Section \ref{ResidueVariety} we directly identify $\HH^1(\overline{C^{*\bullet}_{\X}})$ with the Picard group of the (scheme-theoretic) perfection of $X_0$ (the toric variety with the same fan over the special fiber).  Since Picard groups of schemes arising as inverse limits are easily computed using classical techniques, this gives the result.}
\end{itemize}

In the three steps outlined above we do not work directly with sheaf cohomology groups, but rather the \v{C}ech cohomology groups associated to the cover $\fU$ induced by the fan.  This allows us to do very explicit computations with cochain complexes.  More precisely, we work with quotients of the \v{C}ech complex for $\cO^{+*}_{\cX^\perf}$, which are not necessarily the same as \v{C}ech complexes associated to the quotient sheaves.  We will therefore prove Theorem \ref{main2}} by establishing the following chain of isomorphisms:
\begin{eqnarray}
  \Pic{\X}&\cong&\HH^1({\X},\cO_{{\X}}^*)\label{iso1}\\
  &\cong&\cHH^1(\fU,\cO_{\X}^{*})\label{iso2}\\
  &\cong&\cHH^1(\fU,\cO_{\X}^{+*})\label{iso3}\\
  &\cong&{\HH^1(\overline{C^{*\bullet}_{\X}})}\label{iso4}\\
  &\cong&\Pic (X_0)[p^{-1}]\label{iso5}\\
  &\cong&\Pic (X)[p^{-1}]\label{iso6}.
\end{eqnarray}
 \vspace{0.5cm}
\par We know isomorphism (\ref{iso1}) holds for any locally ringed space.  Proposition \ref{trivialPic} shows that all line bundles on $\cU_\sigma^\perf$ are trivial so that the low degree terms of the \v{C}ech-to-derived functor spectral sequence establish isomorphism (\ref{iso2}).  Isomorphism (\ref{iso6}) is the content of Proposition \ref{toricPic}.  It remains to establish the isomorphisms (\ref{iso3}),(\ref{iso4}) and (\ref{iso5}).
\subsection{Comparison to the Sheaf of Integral Elements}\label{integral}
In this section we establish isomorphism (\ref{iso3}), {completing the extension step of our outline}.  In fact, we prove something slightly more general, {from which the desired isomorphism follows setting $i=1$.}
\begin{Lemma}\label{integralGluing}
  For all $i>0$, we have $\cHH^i(\fU,\cO_{\X}^*)\cong\cHH^i(\fU,\cO_{\X}^{+*})$.
\end{Lemma}
\begin{proof}
  For each $\sigma\in\Sigma$, there is a canonical Gauss norm
   \[|\cdot|:K\la\sigma^{\vee}\cap M[p^{-1}]\ra\longto|K^*|\cup\{0\},\]
  which extends the nonarchimedean absolute value of $K$.  It is given by the rule $|\sum\lambda_m\chi^m| = \sup\{|\lambda_m|\}$.  Restricting to unit groups induces a surjective homomorphism $\cO_{\X}^*(\cU_\sigma^\perf)\to|K^*|$ whose kernel consists of invertible elements with absolute value 1.  But these are precisely the invertible integral elements $\cO_{\X}^{+*}(\cU_\sigma^\perf)$.  As the Gauss norm commutes with localization, we have the following exact sequence of cochain complexes:
  \vspace{0.5cm}
  \[
  \begin{tikzcd}
    0\rar& \prod_i\cO_{\X}^{+*}(\cU^{\perf}_i)\rar\dar&\prod_i\cO_{\X}^{*}(\cU^{\perf}_i)\rar\dar&\prod_i|K^*|\rar\dar&0\\
    0\rar& \prod_{i<j}\cO_{\X}^{+*}(\cU^{\perf}_{ij})\rar\dar&\prod_{i<j}\cO_{\X}^{*}(\cU^{\perf}_{ij})\rar\dar&\prod_{i<j}|K^*|\rar\dar&0\\
    0\rar& \prod_{i<j<k}\cO_{\X}^{+*}(\cU^{\perf}_{ijk})\rar\dar&\prod_{i<j<k}\cO_{\X}^{*}(\cU^{\perf}_{ijk})\rar\dar&\prod_{i<j<k}|K^*|\rar\dar&0.\\
    &\vdots&\vdots&\vdots&
  \end{tikzcd}
  \]
  \vspace{0.5cm}
\par The left and middle sequences are the \v{C}ech sequences for $\cO_{\X}^{+*}$ and $\cO_{\X}^*$ respectively.  Furthermore, the sequence on the right has kernel $|K^*|$ and is otherwise exact (arguing for example as in \cite[Tag 02UW]{Stacks}).  Therefore taking the associated long exact sequence of cohomology groups completes the proof.
\end{proof}

\subsection{{Deformations}}\label{deformation}

In this section we establish isomorphism (\ref{iso4}), {completing the deformation step of our outline}. 

\begin{Def}\label{complexDef}
Let $ C_{\X}^{*\bullet}$ denote the \v{C}ech complex for $\cO_{\X}^{+*}$ associated to $\fU$. Let $\overline{C_{\X}^{*\bullet}}$ denote the reduction mod $1+\cO^{++}_{\cX^\perf}$ of $ C_{\X}^{*\bullet}$. In particular, each $\cO^{+*}_{\cX^\perf}(\cU^\perf_\sigma)$ term of $C_{\X}^{*\bullet}$ is replaced by \[\cO^{+*}_{\cX^\perf}(\cU^\perf_\sigma)/(1+\cO^{++}_{\cX^\perf}(\cU^\perf_\sigma))\cong  (\cO^+_{\X}(\cU^\perf_\sigma)/\cO^{++}_{\X}(\cU^\perf_\sigma))^*.\] See \cite[Lemma 3.8]{DHthesis} and the discussion directly after for more details on this isomorphism.
\end{Def}

The reduction map $ C_{\X}^{*\bullet}\rightarrow \overline{C_{\X}^{*\bullet}}$ of complexes induces a natural map on cohomology

\begin{equation}\label{phiDefined}
  \varphi:\HH^1( C_{\X}^{*\bullet})\longto{\HH^1(\overline{C_{\X}^{*\bullet}})},
\end{equation}

and the main result of this section is the following:

\begin{Prop}\label{defoMain}
  The map $\varphi$ defined in Equation (\ref{phiDefined}) is an isomorphism.
\end{Prop}

As a first step, we establish surjectivity.
\begin{Lemma}\label{sectionOfPhi}
  The map $\varphi$ has a natural section $\psi$.  In particular, $\varphi$ surjects.
\end{Lemma}
\begin{proof}

 {As $\{0\}$ is a face of each $\sigma_i$ there are compatible embeddings \[\cO^+_{\cX^\perf}(\cU_i^\perf)\into K^\circ\left[x_1^{\pm1/p^\infty},\cdots,x_n^{\pm1/p^\infty}\right]\]} reducing to compatible embeddings  \[\cO^+_{\cX^\perf}(\cU_i^\perf)/\cO^{++}_{\cX^\perf}(\cU^\perf_i)\into k\left[x_1^{\pm1/p^\infty},\cdots,x_n^{\pm1/p^\infty}\right].\]
 In particular, we may give \v{C}ech cocycles those coordinates and {observe that} invertible elements must be monomials.  {Denote by $\overline\delta^i$ the differentials of the \v{C}ech complex $\overline{C_{\X}^{*\bullet}}$,} and fix a cocycle
  \[\alpha = (\alpha_{ij}) = (\lambda_{ij}x_1^{m_{ij,1}}\cdots x_n^{m_{ij,n}})\in {\ker\overline\delta^1}.\]
  We first notice that the cocycle of constants $(\lambda_{ij})$ is a coboundary.  Indeed, letting \[x = (\lambda_{1,r},\lambda_{2,r}\cdots,\lambda_{r-1,r},1),\]then the cocycle condition on the $\lambda_{ij}$ implies $\delta^0(x) = (\lambda_{ij})$.  In particular, the class of $\alpha$ in cohomology does not depend on the constant terms of the monomial, so we may take $\alpha_{ij} = x_1^{m_{ij,1}}\cdots x_n^{m_{ij,n}}$.  We therefore define $\psi$ on the level of cocycles by
  \[\psi\left((x_1^{m_{ij,1}}\cdots x_n^{m_{ij,n}})\right)=(x_1^{m_{ij,1}}\cdots x_n^{m_{ij,n}})\] which is certainly a section of $\varphi$.
\end{proof}

To show that this section is in fact an inverse, we construct some intermediary \v{C}ech complexes to interpolate between $C^{*\bullet}_{\X}$ and {$\overline{C^{*\bullet}_{\X}}$}. Let $\varpi\in K$ be a pseudouniformizer.  Then for each positive $d\in\Z[p^{-1}]$ there is a sheaf of principal ideals $(\varpi^d)\hookrightarrow\cO_{\X}^+$. This induces an injection $1+(\varpi^d)\hookrightarrow \cO_{\X}^{+*}$, which gives an injection of \v{C}ech complexes $\check C(\fU, 1+(\varpi^d))\hookrightarrow C^{*\bullet}_{\X}.$ To ease notation, we write $(1+(\varpi^d))^\bullet$ for the complex $\check C(\fU, 1+(\varpi^d))$. 

Let $ C_{\X,d}^{*\bullet}$ denote the reduction mod $(1+(\varpi^d))^\bullet$ of $ C_{\X}^{*\bullet}$. So we have the short exact sequence of complexes 
\[1\rightarrow (1+(\varpi^d))^\bullet \rightarrow C^{* \bullet}_{\X}\rightarrow C^{*\bullet}_{\X,d}\rightarrow 1.\]  
Somewhat more explicitly, we have \[ C_{\X,d}^{*\bullet}=0\rightarrow \displaystyle\prod_i \cO_{\cX^\perf}^{+*}(\cU_i^\perf)/(1+(\varpi^d))\rightarrow \displaystyle\prod_{i<j} \cO_{\cX^\perf}^{+*}(\cU_{ij}^\perf)/(1+(\varpi^d)) \rightarrow \cdots.\] This is a complex of $K^\circ/\varpi^d$-modules. We have natural maps $ C_{\X,d'}^{*\bullet}\rightarrow  C_{\X,d}^{*\bullet}$ whenever $d'>d$.

{We show that the complexes $ C_{\X,d}^{*\bullet}$ interpolate continuously between $\overline{C_{\X}^{*\bullet}}$ and $ C_{\X}^{*\bullet}$. 

\begin{Lemma}[\cite{DHthesis} Lemma 3.10]\label{directLimit}
  $\dlim  C_{\X,d}^{*\bullet} \cong \overline {C_{\X}^{*\bullet}}$ and $\ilim C_{\X,d}^{*\bullet}\cong C_{\X}^{*\bullet}$.
\end{Lemma}
\begin{proof}
To compute the direct limit, we show that
\begin{equation}\label{unionStatement}
  \dlim_{d\in\bZ[p^{-1}]_{>0}}(1+(\varpi^d)) = 1 +\cO_{\X}^{++}.
\end{equation}
We interpret the colimit as a union, and notice that Equation (\ref{unionStatement}) follows if
\[\cO_{\X}^{++} = \bigcup_{d\in\bZ[p^{-1}]_{>0}}(\varpi^d).\]
To see this, we fix an affinoid open $U\subseteq X$, and a topologically nilpotent function $f\in\cO_{\X}^{++}(U)$.  Since $f$ is topologically nilpotent, we can take large $r$ so that $f^{p^r}$ lands in the ideal of $\cO_{\X}^+(U)$ generated by $\varpi$ (as this ideal is an open neighborhood of 0).  We write $f^{p^r} = g\varpi$ for $g\in\cO_{\X}^+(U)$.  Rephrasing we see
\[\left(\frac{f}{\varpi^{1/p^r}}\right)^{p^r}\in\cO_{\X}^+(U).\]
But $\cO_{\X}^+(U)$ is integrally closed in $\cO_{\X}(U)$, so this in turn implies that $f/\varpi^{1/p^r}\in\cO_{\X}^+(U)$, or equivalently, that $f$ is contained in the ideal of $\cO_{\X}^+(U)$ generated by $\varpi^{1/p^r}$.  In particular, $f$ is contained in the union of the ideals generated by the $\varpi^d$ for $d\in\bZ[p^{-1}]_{>0}$.  This completes the verification of Equation (\ref{unionStatement}).  

For the inverse limit, notice that $\ilim C_{\X,d}^\bullet\cong C_{\X}^\bullet$ since $C_{\X}^\bullet$ is $\varpi$-adically complete. Since the unit group functor commutes with inverse limits (indeed, it is right adjoint to the group ring functor), we are done.
\end{proof}
}

{We now will deform the \v{C}ech cohomology groups in question along these quotient complexes.  The following lemma sets up the argument. 

\begin{Lemma}[\cite{DHthesis} Lemma 3.12]\label{l1}
	Let $C_{\X}^\bullet$ be the \v{C}ech complex for $\cO^+_{\cX^\perf}$ associated to $\cU$; let $C_{\X,d}^\bullet$ be the quotient by $(\varpi^d)^\bullet$. For all $d\in\bZ[p^{-1}]_{>0}$ and $i>0$, $\HH^i(C_{\X,d}^\bullet)=0$.
\end{Lemma}
\begin{proof}
	Consider the long exact sequence of \v{C}ech complexes induced by the short exact sequence $0\rightarrow (\varpi^d)^\bullet\rightarrow C_{\X}^\bullet \rightarrow C^\bullet_{\X,d}\rightarrow 0$.  Since $(\varpi^d)\cong\cO_{\X}^+$, Proposition \ref{intAcyclic} tells us that the left and middle complexes are acyclic, so the right one is as well.
\end{proof}
This allows us to show that changing $d$ leaves the cohomology of $C_{\X,d}^{*\bullet}$ unchanged, allowing us to deform.
\begin{Lemma}[\cite{DHthesis} Lemmas 3.13 and 3.14]\label{l2}
For all $d'>d>0$ in $\bZ[p^{-1}]_{>0}$ and $i>0$, the natural map
\[\HH^i(C_{\X,d'}^{*\bullet})\longto\HH^i(C_{\X,d}^{*\bullet})\]
is an isomorphism.  When $i=0$, the map\[\HH^0(C_{\X,d'}^{*\bullet})\longto\HH^0(C_{\X,d}^{*\bullet})\] is surjective.
\end{Lemma}
}
\begin{proof}
{
We first consider the case where $d'=2d$.  We have the following diagram whose rows are exact sequences of \v{C}ech complexes,
}
\[
\begin{tikzcd}
  0\rar & (1+(\varpi^{2d}))^\bullet\rar\ar[hookrightarrow]{d} & C_{\X}^{*\bullet} \ar[equal]{d} \rar & C^{*\bullet}_{\X,2d}\rar\ar[twoheadrightarrow]{d} & 0\\
  0\rar & (1 + (\varpi^d))^\bullet\rar & C_{\X}^{*\bullet}\rar & C^{*\bullet}_{\X,d}\rar & 0
\end{tikzcd}
\]
{
The snake lemma exhibits the exact sequence
\begin{equation}\label{*5}
	1\longto \frac{(1+(\varpi^d))^\bullet}{(1+(\varpi^{2d}))^\bullet}\longto C^{*\bullet}_{\X,2d} \longto C^{*\bullet}_{\X,d} \longto1.
\end{equation}
We claim that the quotient of complexes $\frac{(1+(\varpi^d))^\bullet}{(1+(\varpi^{2d}))^\bullet}$ is isomorphic to $C^\bullet_{\X,d}$. To see this, note that for any open subset $U$ we have:
\[\frac{(1+(\varpi^d))^\bullet}{(1+(\varpi^{2d}))^\bullet}(U)\cong\frac{1+\varpi^d\cO_{\X}^+(U)}{1+\varpi^{2d}\cO_{\X}^+(U)}\cong 1+\frac{\varpi^d\cO_{\X}^+(U)}{\varpi^{2d}\cO_{\X}^+(U)}.\]
We can then construct a map to the right hand side from $\cO_{\X}^+(U)$ by the rule $a\mapsto1+a\varpi^d$.  The map is a homomorphism because $\varpi^d/\varpi^{2d}$ squares to zero, it is clearly surjective, and the kernel is precisely the ideal generated by $\varpi^d$ so that it induces an isomorphism from $C^\bullet_{\X,d}$.  In particular, by Lemma \ref{l1}, $1+\varpi^d/\varpi^{2d}$ has no higher \v{C}ech cohomology.  Therefore the case of the result for $d'=2d$ follows from the long exact sequence on \v Cech cohomology associated Sequence (\ref{*5}).

For general $d'>d>0$, the previous paragraph allows us to replace $d$ with $2^ld$, so we may assume $d<d'<2d$.  Then when $i>0$, we have the following commutative diagram, where the top and bottom rows are isomorphisms again by the previous paragraph.
}
\[
\begin{tikzcd}
  {} & \HH^i(C^{*\bullet}_{\X,2d})\ar["\gamma"]{dr}\ar["\sim"]{rr} & {} & \HH^i(C^{*\bullet}_{\X,d})\\
  \HH^i(C^{*\bullet}_{\X,2d'})\ar{ur}\ar["\sim"]{rr} & {} & \HH^i(C^{*\bullet}_{\X,d'})\ar{ur} & {}
\end{tikzcd}
\]
{
In particular, $\gamma$ is injective and surjective, hence an isomorphism. If $i=0$, the bottom row is surjective, so $\gamma$ must be surjective as desired.
}
\end{proof}

\begin{Lemma}\label{limitCommutes}
The map $\HH^1(C^{*\bullet}_{\X})\rightarrow \ilim \HH^1(C^{*\bullet}_{\X,d})$ is an isomorphism.
\end{Lemma}

\begin{proof}

Form the inverse system $ C_{\X}^*=( C_{\X,d}^{*\bullet})_{d\in \bN}$. This is an object of the derived category of $\operatorname{Mod}(\bN,(K^\circ/\varpi^d))$, where $\operatorname{Mod}(\bN,(K^\circ/\varpi^d))$ is the category of inverse systems $(M_d)$ of abelian groups such that each $M_d$ is a $K^\circ/\varpi^d$-modules and the transition maps $M_{d+1}\rightarrow M_d$ is a map of $K^\circ/\varpi^{d+1}$-modules.

By \cite[Tag 0CQE]{Stacks}, we have a short exact sequence of $K^\circ$-modules:

\begin{equation}\label{eq:ses}
0\rightarrow R^1\varprojlim_d \HH^0(C_{\X,d}^{*\bullet})\rightarrow \HH^1(R\lim  C_{\X}^*)\rightarrow \varprojlim_d \HH^1(C_{\X,d}^{*\bullet})\rightarrow 0. 
\end{equation}

We first check that $R\lim  C_{\X}^* \cong  C_{\X}^{*\bullet},$ showing the right arrow in the sequence the map we care about.  

Consider the diagram

\begin{tikzcd}
0 \arrow[r] & \displaystyle\prod_i \cO_{\cX^\perf}^{+*}(\cU_i^\perf) \arrow[r]\arrow[d] & \displaystyle\prod_{i<j} \cO_{\cX^\perf}^{+*}(\cU_{ij}^\perf) \arrow[r] \arrow[d]& \cdots \\
 & \vdots \arrow[d]& \vdots \arrow[d]& \cdots \\
0 \arrow[r] & \displaystyle\prod_i \cO_{\cX^\perf}^{+*}(\cU_i^\perf)/(1+(\varpi^2)) \arrow[d]\arrow[r]& \displaystyle\prod_{i<j} \cO_{\cX^\perf}^{+*}(\cU_{ij}^\perf)/(1+(\varpi^2)) \arrow[d] \arrow[r] & \cdots \\
0 \arrow[r] & \displaystyle\prod_i \cO_{\cX^\perf}^{+*}(\cU_i^\perf)/(1+(\varpi))\arrow[r] & \displaystyle\prod_{i<j} \cO_{\cX^\perf}^{+*}(\cU_{ij}^\perf)/(1+(\varpi)) \arrow[r] & \cdots .
\end{tikzcd}

Here the top row is $ C_{\X}^{*\bullet}$ and everything below it is $ C_{\X}^*$. By Lemma \ref{directLimit} we see that each term in the top row is the inverse limit of the modules in that column. 

All vertical maps are quotients, so the columns satisfy the Mittag-Leffler condition. Let $ C_{\X,d}^{*,p}$ be the $p$th term of the complex $ C_{\X,d}^{*\bullet}$. Then by \cite[Tag 091D(3)]{Stacks}, $R^1\varprojlim_{d}  C_{\X,d}^{*,p}=0$ for all $p$. Using this, \cite[Tag 091D(5)]{Stacks} implies that $R\lim  C$ is represented by the complex whose term in degree $p$ is $\varprojlim_{d}  C_{\X,d}^{*,p}$, which is exactly the top row of the diagram. 

To show that the right map in Equation \ref{eq:ses} is an isomorphism, we must show that $R^1\varprojlim_d \HH^0(C_{\X,d}^{*\bullet})$ vanishes. By Lemma \ref{l2}, the inverse system $(\HH^0(C_{\X,d}^{*\bullet}))_{d\in \bN}$ has surjective transition maps. It is therefore Mittag-Leffler, so its $R^1\lim$ vanishes as desired.

\end{proof}

We are now ready to assemble all this work into a proof of Proposition \ref{defoMain}. 

\begin{proof}[Proof of Proposition \ref{defoMain}]
We have shown that $\varphi$ decomposes into the following composition of isomorphisms,
\begin{eqnarray}
  \HH^1(C_{\X}^{*\bullet})&\cong& {\HH^1(\ilim C^{*\bullet}_{\X,d})}\label{ea1}\\
  &\cong&{\ilim\HH^1(C^{*\bullet}_{\X,d})}\label{ea2}\\
  &\cong& {\HH^1(C^{*\bullet}_{\X,d})}\label{ea3}\\
  &\cong&{\dlim \HH^1(C^{*\bullet}_{\X,d})}\label{ea4}\\
  &\cong& {\HH^1(\dlim C^{*\bullet}_{\X,d})}\label{ea5}\\
  &\cong&{\HH^1(\overline{C_{\X}^{*\bullet}})}\label{ea6}.
\end{eqnarray}
{That lines (\ref{ea1}) and (\ref{ea6}) are isomorphisms is Lemma \ref{directLimit}. That line (\ref{ea2}) is an isomorphism is Lemma \ref{limitCommutes}. Lines (\ref{ea3}) and (\ref{ea4}) are isomorphisms because they are inverse and direct limits of a system of isomophisms due to Lemma \ref{l2}. Line (\ref{ea5}) is an isomorphims because filtered colimits of abelian groups are exact.

}

\end{proof}

As $\cHH^1(\fU,\cO^{+*}_{\X})= \HH^1(C^{*\bullet}_{\X})$ by definition, we have now established isomorphism (\ref{iso4}).

\subsection{Comparison to the Perfection of the Toric Variety over the Residue}\label{ResidueVariety}
To finish we must establish isomorphism \eqref{iso5}, between {$\HH^1(\overline{C^{*\bullet}_{\X}})$} and $\Pic(X_0)\left[p^{-1}\right]$.  To do so, we establish a bit of notation.  Let $U_i = U_{\sigma_i,k}$, and define $\fU_0 = \{U_i\to X_0\}$ the standard cover of the toric variety $X_0$ induced by $\Sigma$.  Let $X_0^\sperf$ be the (scheme theoretic) perfect closure of $X_0$ (that is, the inverse limit along Frobenius in the category of schemes). Then $\fU_0^\sperf = \{U_i^\sperf\to X_0^\sperf\}$ is an open cover of $X_0^\sperf$ as perfect closures commute with base change.
\begin{Lemma}
  There is a natural isomorphism {$\HH^1(\overline{C^{*\bullet}_{\X}})\cong\cHH^*(\fU_0^\sperf,\cO_{X^\sperf_0}^*)$}.
\end{Lemma}
\begin{proof}
  Since $k$ is perfect, there is a natural identification
  \[ k\left[\sigma^\vee\cap M[p^{-1}]\right] = \HH^0(U_\sigma^\sperf,\cO_{X^\sperf_0}).\]
  In particular the unit groups and localization maps are identified, so that this passes to a natural identification of complexes
  \[\overline{C^{*\bullet}_{\X}} \cong \check C(\fU_0^\sperf,\cO_{X^\sperf_0}^*).\]
  Passing to cohomology completes the proof.
\end{proof}
By \cite[Theorem 2.1]{Gub}, every line bundle on $X_0$ trivializes on $\fU_0$, so that the \v{C}ech-to-derived functor spectral sequence gives an isomorphism $\cHH^1(\fU_0,\cO_{X_0}^*)\cong\Pic X_0$.  Furthermore, we can pass \cite[Theorem 2.1]{Gub} to the colimit along Frobenius (applying, for example \cite[Tag 0B8W]{Stacks}), so that every line bundle on $X_0^\sperf$ trivializes on $\fU_0^\sperf$.  Again by the \v{C}ech-to-derived functor spectral sequence we see that $\cHH^1(\fU_0^\sperf,\cO_{X^\sperf_0}^*)\cong\Pic X_0^\sperf$.
\begin{Prop}\label{perfectPic}
  \textcolor{black}{Suppose $S$ is a quasicompact quasiseparated scheme over a field of characteristic $p$.  Let $S^\sperf$ be its perfect closure. Then $\Pic(S^\sperf)\cong\Pic(S)[p^{-1}]$.}
\end{Prop}
\begin{proof}
  \textcolor{black}{As $S^\sperf$ is perfect, the Frobenius map on $\cO^*_{S^\sperf}$ is an isomorphism.  Passing to the first cohomology group shows that the $p$-th power map on $\Pic(S^\sperf)$ is an isomorphism.  This implies that the pullback map $\Pic(S)\to\Pic(S^\sperf)$ factors through $\Pic(S)[p^{-1}]$, so we have,}
  \[\textcolor{black}{\Pic(S)\to\Pic(S)[p^{-1}]\to\Pic(S^\sperf).}\]
  \textcolor{black}{We must show the righthand map is an isomorphism.  We remind ourselves that $S^\sperf = \ilim S$ along Frobenius.  As Frobenius is an affine morphism and $S$ is quasicompact and quasiseparated, every line bundle on $S^\sperf$ is the pullback of a line bundle from one of the factors \cite[Tag 0B8W]{Stacks}.  This proves surjectivitiy.}

	\textcolor{black}{To show injectivity it is equivalent to show that the kernel of $\Pic(S)\to\Pic(S^\sperf)$ is $p$-power torsion.  Fix some $\cL$ in the kernel, and represent it by a cover $\{V_i\to S\}$ and gluing functions $f_{jk}\in\cO_S(V_{jk})^*\subseteq\cO_{S^\sperf}(V_{jk}^\sperf)^*$.  Then the fact that $\cL$ pulls back to a trivial bundle means that $f_{jk} = f_j/f_k$ for various $f_i\in\cO_{S^\sperf}(V_i^\sperf)^*$.  There is some large $N$ such that each for each $i$ we have $f_i^{p^N}\in\cO_S(V_i)^*$, so that $f_{jk}^{p^N} = f_j^{p^N}/f_k^{p^N}$.  In particular the cocycle $(f_{jk}^{p^N})$ is a coboundary and therefore its gluing data is trivial, whence $\cL^{\otimes p^N}\cong\cO_{S^{\sperf}}$.}
\end{proof}
As $X_0$ is proper, we have established that $\Pic(X_0^\sperf)\cong\Pic(X_0)[p^{-1}]$, which is isomorphism \ref{iso5}.  This was the final link in the proof of Theorem \ref{main2}.
\begin{Remark}[A note on canonicity]\label{canon}
  Notice that the multiplication by $p$ map is an isomorphism on $\Pic({\X})$.  Once we fix a projection map $\pi_0:{\X}\to X$ the isomorphism from Theorem \ref{main2} becomes canonical in the following sense.  The construction of the perfectoid cover produces the following commutative diagram where $\varphi = \varphi_p$ is the $p$th power map of Definition \ref{powerMap}
  \vspace{0.25cm}
  \[
  \begin{tikzcd}
    {\X}\rar\ar[rr,bend left = 30,swap,"\pi_{k+1}"]\ar[rrr,bend left = 30,"\pi_k"]\ar[rrrrr,bend left=30,"\pi_0"] & \cdots\rar & X\ar[r,"\varphi"] & X\ar[r] & \cdots\rar & X.
  \end{tikzcd}
  \]
  \vspace{0.25cm}\\
  Passing to Picard groups this induces by universal property a canonical homomorphism
	\[\dlim_{\varphi^*}\Pic(X)\to\Pic({\X}).\]
	By Corollary \ref{powerPullback}, $\varphi^*$ is multiplication by $p$ so that the source is canonically isomorphic to $\Pic(X)[p^{-1}]$.  Composing gives a canonical homomorphism
	\[\Pic(X)[p^{-1}]\to\Pic({\X})\]
	which is an isomorphism by Theorem \ref{main2}.  In particular, identifying $\cL\in\Pic({\X})$ with a formal $p$th tensor root $\cM^{1/p^k}$ of some $\cM\in\Pic(X)$ in turn identifies $\cL$ with the pullback of $\cM$ from the $k$th level of the tower: $\cL\cong\pi_k^*\cM$.
\end{Remark}
\section{Cohomology of Line Bundles}
We will conclude with a computation of the cohomology of line bundles on the perfectoid cover of a toric variety.  The standard setup will be the following.
\begin{Setup}\label{cohomologySetup}
	Let $\Sigma$ be a complete fan, $K$ a perfectoid field, and $X = X_{\Sigma,K}$ the associated toric variety with perfectoid cover ${\X}\to X$. Assume that the Picard group of $\cX_\Sigma$ trivializes on the affinoid cover induced by $\Sigma$ (so that Theorem \ref{main2} applies).  As in Remark \ref{canon}, fix $\cL\in\Pic{\X}$ and a line bundle $\cM\in\Pic X$ such that:
	\[\cL\cong\cM^{1/p^k}\cong\pi_k^*\cM,\]
	where $\pi_k$ is projection onto the $k$-th factor of the inverse limit.
\end{Setup}
\begin{Th}\label{cohomologyTheorem}
	In the situation of Setup \ref{cohomologySetup}, for all $i\ge0$, there is a canonical homomorphism
	\[\gamma:\operatorname{colim}_n\HH^i(X,\cM^{p^{n-k}}){ }\longto\HH^i({\X},\cL).\]
	Furthermore, one can endow the source and target with with the structure of topological $K$-vector spaces in a way such that the target is the completion of the source, and $\gamma$ is the canonical inclusion.  In particular, with this topology fixed, we have:
	\[\left(\operatorname{colim}_n\HH^i(X,\cM^{p^{n-k}})\right)^\wedge{ }\cong\HH^i({\X},\cL).\]
\end{Th}
\begin{Remark}\label{k=0}
	Without loss of generality we may assume that the identification $\cL = \cM^{1/p^k}$ can be made with $k=0$.  This is because the $p$-power map on $\varphi:\X\to\X$ (Definition \ref{powerMap}) is an isomorphism, as it is determined affinoid locally by multiplication by $p$ on the free abelian group $M[1/p]$.  The isomorphism $\varphi^*\cL\cong\cL^p$, (Corollary \ref{powerPullback}) is adjoint to an isomorphism $\cL\cong\varphi_*\cL^p$, which as $\varphi$ is an isomorphism passes to a $K$-linear isomorphism of the cohomology of $\cL$ and $\varphi_*\cL^p$.  This composes to a canonical isomorphism on the cohomology of $\cL$ and $\cL^p$ applying Lemma \ref{cohomologyPushforward} below.
\end{Remark}
\subsection{Construction of the Comparison Map}
We first construct the canonical map.  We need the following immediate consequence of the Leray spectral sequence.
\begin{Lemma}\label{cohomologyPushforward}
	Let $f:Y\to Z$ be a map of topological spaces, and $\cF$ a sheaf of abelian groups on $Y$.  For all $i$ there are natural maps $\HH^i(Z,f_*\cF)\to\HH^i(Y,\cF)$, which are isomorphisms if $f$ is.
\end{Lemma}
\begin{proof}
	This is an immediate consequence of the filtration of the $E_\infty$ page from the Leray spectral sequence
	\[E^{p,q}_2:\hspace{20pt}\HH^p(Z,R^qf_*\cF)\Longrightarrow \HH^{p+q}(Y,\cF).\]
\end{proof}
We will now build the map from Theorem \ref{cohomologyTheorem}.  We assume the conditions of Setup \ref{cohomologySetup}, identifying $\cL\in\Pic\X$ with $\cM^{1/p^k}$, for $\cM\in\Pic\cX$, By Remark \ref{k=0} we may assume $k=0$.
\begin{Prop}\label{gammaHat}
	There is a canonical homomorphism
	\[\gamma:\operatorname{colim}_n\HH^i(X,\cM^{p^n})\longto\HH^i({\X},\cL).\]
\end{Prop}
\begin{proof}
	Let $\varphi$ be the $p$th power map on $X$. We know by Corollary \ref{powerPullback} that there is an isomorphism $\varphi^*\cM\cong\cM^p$,  which is adjoint to a map $\cM\to\varphi_*\cM^p$.  Passing to cohomology and composing with the map from Lemma \ref{cohomologyPushforward} gives a homomorphism $\rho:\HH^i(X,\cM)\to\HH^i(X,\cM^p).$  Arguing similarly, for each  $m>0$, the isomorphisms $\cL\cong\pi_{m}^*\cM^{p^m}$ induce canonical maps $\gamma_{m}:\HH^i(X,\cM^{p^m})\to\HH^i({\X},\cL)$, and these fit compatibly in the following diagram.
	\begin{equation}\label{cohomologyDiagram}
	\begin{tikzcd}
	  \HH^i(X,\cM)\ar[d,swap,"\rho"]\ar[dddrr,"\gamma_0"]&&\\
	  \HH^i(X,\cM^{p})\dar\ar[ddrr,swap,"\gamma_{1}"]&&\\
	  \HH^i(X,\cM^{p^2})\dar\ar[drr,swap,"\gamma_{2}"]&&\\
	  \vdots\dar&&\HH^i({\X},\cL)\\
	  \operatorname{colim}\HH^i(X,\cM^{p^n}).\ar[urr,dotted,"\exists!\gamma"]
	\end{tikzcd}
	\end{equation}
\end{proof}
\subsection{Topologizing the Cohomology Groups}
We will use \v{C}ech cohomology with respect to the usual cover to endow the source and target of $\gamma$ with topologies.  We let $\{\sigma_1,\cdots,\sigma_r\}$ be the maximal cones of $\Sigma$, and consider the covers $\fU = \{U_{\sigma_i}\to X\}$ and $\fU^\perf = \{\cU^\perf_{\sigma_i}\to{\X}\}$.  We first record that \v{C}ech cohomology is effective with respect to these covers.
\begin{Lemma}\label{lbCechEffective}
Consider the situation of Setup \ref{cohomologySetup}, and let $\fU$ and $\fU^\perf$ be the standard covers of $X$ and $\X$ respectively.  Then the natural maps:
\[\cHH^i(\fU^\perf,\cL)\longtoo{\sim}\HH^i(\cX^\perf,\cL)\hspace{10pt}\text{ and }\hspace{10pt}\cHH^i(\fU,\cM^{p^n})\longtoo{\sim}\HH^i(X,M^{p^n}),\]
are isomorphisms.
\end{Lemma}
\begin{proof}
The $U_{\sigma_i}$ and their intersections are affinoid and  and the same can be said for the $\cU^\perf_{\sigma_i}$.  Since locally free sheaves on affinoid adic spaces are acyclic, the \v{C}ech-to-derived functor spectral sequence gives isomorphisms \ref{cIso1} and \ref{cIso5}.
\end{proof}
Now we can explicitly write down the maps $\rho$ and $\gamma_i$ from the proof of Proposition \ref{gammaHat} on \v{C}ech cocycles.
\begin{Lemma}\label{explicit}
	For all maximal cones $\sigma\in\Sigma$, the endomorphism of $K[\sigma^{\vee}\cap M]$ prescribed by $\chi^m\mapsto\chi^{pm}$ induces maps $\cO_{U_{\sigma_i}}\to\cO_{U_{\sigma_i}}$. These maps glue to the map $\cM\to\varphi_*\cM^p$.
\end{Lemma}
\begin{proof}
  We know that $\cM$ is given by a cocycle $(\chi^{m_{ij}})\in\cHH^1(\fU,\cO_X^*)$, (for some $m_{ij}\in\sigma_{ij}^\vee\cap M$), and $\cM^p$ is given by $(\chi^{pm_{ij}})$.  Therefore it suffices to show that the given map commutes with the gluing data on $U_{ij}$.  But one can easily check that
  \[
  \begin{tikzcd}
    \cO_{U_{ij}}\ar[d,swap,"\cdot\chi^{m_{ij}}"]\ar[rr,"\chi^m\mapsto\chi^{pm}"]&&\cO_{U_{ij}}\ar[d,"\cdot\chi^{pm_{ij}}"]\\
    \cO_{U_{ji}}\ar[rr,"\chi^m\mapsto\chi^{pm}"]&&\cO_{U_{ji}}
  \end{tikzcd}
  \]
  commutes, so we are done.
\end{proof}
This has the following immediate consequence:
\begin{Lemma}\label{explicit2}
	For any $\sigma\in\Sigma$, the map $\HH^0(U_\sigma,\cM)\to \HH^0(U_\sigma,\cM^p)$ can be identified with the inclusion
	\[K[\sigma^\vee\cap M]\into K\left[\sigma^\vee\cap p^{-1}M\right],\]
	in such a way that it is compatible with restrictions to the faces of $\sigma$.  In particular, one can identify:
	\[\colim \HH^0(U_\sigma,M^{p^n})\cong K[\sigma^\vee\cap M[p^{-1}]],\]
	compatibly with restricting to faces in $\sigma$.
\end{Lemma}
\begin{proof}
	This follows immediately from Lemma \ref{explicit} relabelling $ pM \subseteq M$ as $M\subseteq p^{-1}M$
\end{proof}
Diagram \ref{cohomologyDiagram} from the proof of Proposition \ref{gammaHat} can be exhibited as the induced map on the cohomology of the following composition of \v{C}ech complexes:
  \vspace{0.25cm}
	\[\check C(\fU,\cM)\to \check C(\fU,\cM^p)\to\cdots\to\operatorname{colim}\check C(\fU,\cM^{p^n})\too{\eta}\check C(\fU^\perf,\cL)\]
  \vspace{0.25cm}\\
	To simplify notation, we give the source and target of $\eta$ the names $C^*$ and $D^*$ respectively.  With Lemma \ref{explicit2} in mind, the $\eta$ can be identified with the following inclusion (now arranged vertically):
	\begin{equation}\label{etaDiagram}
	\begin{tikzcd}
	  C^*:=\cdots\rar & \prod_{j_{0...i}}K\left[\sigma_{j_{0...i}}^\vee\cap M[p^{-1}]\right]\rar\dar&\prod_{j_{0...i+1}}K\left[\sigma_{j_{0...i+1}}^\vee \cap M[p^{-1}]\right]\rar\dar&\cdots\\
	  D^*:=\cdots\rar&\prod_{j_{0...i}}K\left\la \sigma_{j_{0...i}}^\vee\cap M[p^{-1}]\right\ra\rar&\prod_{j_{0...i+1}}K\left\la\sigma_{j_{0...i+1}}^\vee\cap M[p^{-1}]\right\ra\rar&\cdots.\\
	\end{tikzcd}
	\end{equation}
This looks like the inclusion of a complex of topological $K$-vector spaces into its termwise completion.  Let's make that precise.
\begin{Prop}\label{summary}
Let $C^*$ and $D^*$ be as above.  Then $C^*,D^*$ and their cohomologies $\HH^i(C^*), \HH^i(D^*)$ can be equipped with topologies under which:
\begin{enumerate}
\item $\eta:C^*\into D^*$ can be identified with the map from a complex of topological $K$-vector spaces to its (termwise) completion.
\item For each index $i$, the map $\gamma:\HH^i(C^*)\to\HH^i(D^*)$ can be identified with the map from a topological $K$-vector space to its completion.
\end{enumerate}
\end{Prop}
\begin{proof}
Because the differentials of the complexes above are generated by alternating sums of monic monomials, they take take polynomials (resp. power series) with integral to ones with integral coefficients.  In particular, Diagram \ref{etaDiagram} restricts to a morphism $\eta^\circ$ of complexes:
\[
	\begin{tikzcd}
	  C^{*,\circ}: \cdots\rar & \prod_{j_{0...i}}K^\circ\left[\sigma_{j_{0...i}}^\vee\cap M[p^{-1}]\right]\rar\dar&\prod_{j_{0...i+1}}K^\circ\left[\sigma_{j_{0...i+1}}^\vee \cap M[p^{-1}]\right]\rar\dar&\cdots\\
	  D^{*,\circ}: \cdots\rar&\prod_{j_{0...i}}K^\circ\left\la \sigma_{j_{0...i}}^\vee\cap M[p^{-1}]\right\ra\rar&\prod_{j_{0...i+1}}K^\circ\left\la\sigma_{j_{0...i+1}}^\vee\cap M[p^{-1}]\right\ra\rar&\cdots.\\
	\end{tikzcd}
\]
Giving the source and target the $\varpi$-adic topology, we see that $\eta^\circ$ is the inclusion of a complex of topological $K^\circ$-modules into its $\varpi$-adic completion.  Now, for each $i$, we have that:
\[C^i = C^{i,\circ}\otimes_{K^\circ}K = C^{i,\circ}[1/\varpi],\]
and we can therefore give $C^i$ the topology making $C^{i,\circ}$ open and bounded (that is, the topology generated by the subsets $\varpi^dC^{i,\circ}$ as $d$ varies), and we do similarly for each $D^k$.  With this topology, then it is clear $\eta$ is the inclusion of a topological $K$-vector space into its completion (as it is true in a neighborhood of 0), proving part (1).

Now on to cohomology; the groups $\HH^i(C^{*,\circ})$ are naturally $K^\circ$-modules and can therefore be endowed with the $\varpi$-adic topology, and similarly for the $\HH^i(D^{*\circ})$.  Since $K$ is flat over $K^\circ$, we have isomorphisms
\[\HH^i(C^*) \cong\HH^i(C^{*,\circ}\otimes_{K^\circ} K) \cong \HH^i(C^{*,\circ})\otimes_{K^\circ} K \cong \HH^i(C^{*,\circ})[1/\varpi].\]
Then we can give $\HH^i(C^*)$ the topology making the image of $\HH^i(C^{*,\circ})$ open and bounded, that is, the topology induced by the images of $\varpi^d\HH^i(C^{*,\circ})$ in $\HH^i(C^*)$ as $d$ varies.  We do similarly for $\HH^i(D^*)$.  With these topologies in mind, $\gamma$ is obtained from the inclusion of $C^*$ into its completion $D^*$ by passage to cohomology.  Part (2) then follows from the more general Lemma \ref{completionExact} below.
\end{proof}
To conclude the proof of we'd like to argue$-$as in the proof of Proposition \ref{intAcyclic}$-$that completion of a sequence of linearly topologized modules commutes with cohomology, (and in fact, this argument immediately implies $\HH^i(C^{*,\circ})^\wedge\cong\HH^i(D^{*,\circ})$).  As the modules in question are not linearly topologized, we need the following lemma extending the commutation of cohomology and completion to generic fibers of linearly topologized $K^\circ$-modules.
\begin{Lemma}\label{completionExact}
  Let $M^* = M^0\to\M^1\to\cdots\to M^n$ be a sequence of linearly topologized $K^\circ$-modules with countable systems of fundamental neighborhoods of 0, and for each $i$ give $M^i\otimes K$ the topology making the image of $M^i$ open and bounded.  Then:
  \[\HH^i(M^*\otimes K)^\wedge{}\cong\HH^i\left((M^*\otimes K)^\wedge{}\right).\]
\end{Lemma}
\begin{proof}
  We collect a few ingredients.  First we use that if $N$ is a linearly topologized $K^\circ$-module, then
  \[\hat N\otimes K\cong(N\otimes K)^\wedge{}.\]
  This is rather immediate, as $N$ is open in $N\otimes K$ and the  basis of 0 given by the $\varpi^n N$ is a basis for both the topology of $N$ and $N\otimes K$.  We will also use that $\cdot\otimes K$ is an exact functor, and that $\HH^i(\hat M^*)\cong\HH^i(M^*)^\wedge{}$ since the $M^i$ are linearly topologized with countable neighborhood bases of 0 \cite[Tag 0AS0]{Stacks}.  Putting all this together gives the following chain of isomorphisms which prove the result,
  \begin{eqnarray*}
    \left(\HH^i(M^*\otimes K)\right)^\wedge{}&\cong&\left(\HH^i(M^*)\otimes K\right)^\wedge{}\\
    &\cong&\HH^i(M^*)^\wedge{}\otimes K\\
    &\cong&\HH^i(\hat M^*)\otimes K\\
    &\cong&\HH^i(\hat M^*\otimes K)\\
    &\cong&\HH^i\left((M^*\otimes K)^\wedge{}\right).
  \end{eqnarray*}
\end{proof}
We can now string everything together to prove the main theorem of this section.
\begin{proof}[Proof of Theorem \ref{cohomologyTheorem}]
	We use the $\wedge$ symbol to represent completion with respect to the topologies introduced in the previous paragraph.  The map from Proposition \ref{gammaHat} can be identified with the following chain of isomorphisms.
	\begin{eqnarray}
		\HH^i({\X},\cL)&\cong&\HH^i(\check C(\fU^\perf,\cL))\label{cIso1}\\
		&\cong&\HH^i\left(\operatorname{colim}\left(\check C\left(\fU,\cM^{p^n}\right)\right)^\wedge{}\right)\label{cIso2}\\
		&\cong&\HH^i\left(\operatorname{colim}\left(\check C\left(\fU,\cM^{p^n}\right)\right)\right)^\wedge{}\label{cIso3}\\
		&\cong&\left(\operatorname{colim}\HH^i\left(\check C\left(\fU,\cM^{p^n}\right)\right)\right)^\wedge{}\label{cIso4}\\
		&\cong&\left(\operatorname{colim}\HH^i\left(X,\cM^{p^n}\right)\right)^\wedge{}\label{cIso5}
	\end{eqnarray}
	To conclude we verify that they are all in fact isomorphisms.  Isomorphisms \ref{cIso1} and \ref{cIso5} are Lemma \ref{lbCechEffective} and isomorphism \ref{cIso4} follows from the fact that cohomology commutes with directed colimits of abelian groups (\cite[Tag 00DB]{Stacks}).  Isomorphism \ref{cIso2} is Proposition \ref{summary}(1) and isomorphism \ref{cIso3} is Proposition \ref{summary}(2), which completes the proof.
\end{proof}

\subsection{Vanishing Theorems}
The cohomology of line bundles on toric varieties has been extensively studied, for example in \cite{Demazure},\cite{BB},\cite[9.1-9.4]{CLS}.  Therefore, given a line bundle $\cL$ on ${\X}$, we may first use Theorem \ref{main2} to identify it with $\cM^{1/p^k}$ for a line bundle $\cM$ on $X$, and then use known results about the cohomology of $\cM$ together with Theorem \ref{cohomologyTheorem} to explicitly compute the cohomology of $\cL$.  We will use this philosophy to promote Demazure and Batyrev-Borisov vanishing theorems for toric varieties to the perfectoid setting.  These theorems concern globally generated line bundles, so to deduce them from the classical theorems about toric varieties using Theorem \ref{cohomologyTheorem}, we need the correspondance between line bundles of Theorem \ref{main2} to preserve the property of being globally generated.
\begin{Prop}\label{globallyGenerated}
  Let $\cL\in\Pic\X$ and identify it with $\cM^{1/p^k}$ for $\cM\in\Pic\cX$.  Then $\cL$ is a globally generated $\cO_{\X}$-module if and only if $\cM^{p^t}$ is a globally generated $\cO_{\cX}$-module for $t>>0$.
\end{Prop}
\begin{proof}
	As in Remark \ref{k=0} we may assume $k=0$.  If $\cM^{p^t}$ is globally generated, then we can pick a surjection $\rho:\bigoplus\cO_X\onto\cM^{p^t}$.  Letting $\varphi$ be the $p$-power map then $\pi_0 = \varphi^t\circ\pi_t$.  Therefore (applying Lemma \ref{powerPullback})
	\[\cL \cong \pi_0^*\cM = \pi_t^*(\varphi^t)^*\cM\cong\pi_t^*\cM^{p^t}.\]
	As pulling back by $\pi_t$ is right exact, $\pi_t^*\rho:\bigoplus\cO_{{\X}}\onto\cL$ gives the desired surjection.

  Conversely, suppose that $\cL$ is globally generated.  As $\X$ is compact, $\cL$ may be generated by a finite set of sections $s_1,\cdots,s_n$.  This proof will have 2 parts.  First we will show that the $s_j$ can without loss of generality be assumed to come from global sections of $\cM^{p^t}$ for some $t$ (under the identification of Theorem \ref{cohomologyTheorem} on zeroth cohomology).  We will then show that (perhaps increasing $t$) these $s_j$ generate $\cM^{p^t}$.

  Trivialize $\cL$ over the $\cU^\perf_i$, and consider the $s_j$ as elements of
  \[\cO_{\X}(\cU^\perf_i) = K\left\la\sigma_i^\vee\cap M[p^{-1}]\right\ra.\]
  As $\cU^\perf_i$ is affinoid, the $s_j$ generate the unit ideal of $\cO_\cX(\cU^\perf_i)$.  Therefore there are $a_1,\cdots,a_n\in\cO_{\X}(\cU^\perf_i)$ such that $a_1s_1+\cdots+a_ns_n = 1$.  Applying Theorem \ref{cohomologyTheorem} on zeroth cohomology, there are global sections $\tilde s_j\in\Gamma(\cX,\cM^{p^t})$ which are arbitrarily close to the $s_j$ (increasing $t$ as necessary).  Since there are finitely many $\cU^\perf_i$, we may choose them so that on each $\cU^\perf_i$:
  \[||(a_1s_1+\cdots+a_ns_n)-(a_1\tilde s_1+\cdots+a_n\tilde s_n)||<1.\]
  That is, $\xi := 1-(a_1\tilde s_1+\cdots+a_n\tilde s_n)$ is topologically nilpotent.  Therefore the geometric series for $(1-\xi)^{-1}$ coverges so that $(1-\xi)$ is a unit.  As $1-\xi$ is in the ideal generated by the $\tilde s_j$, they generate the unit ideal.  In particular, we see that the $\tilde s_j$ generate $\cL$.  This completes the first step.

  The second step applies a similar argument, but now to the coefficients.  We use that trivializations of $\cL$ over the $\cU_i^\perf$ and of $\cM$ over $\cU_i$ are compatible, so that we can identify the map $\cM^{p^t}(\cU_i)\to\cL(\cU_i^\perf)$ with the inclusion $K\left\la\sigma_i^\vee\cap p^{-t}M\right\ra\into K\left\la\sigma_i^\vee\cap M[p^{-1}]\right\ra$, and the $\tilde s_j$ can be considered as elements of the target.  We know they generate the unit ideal in $K\left\la\sigma_i^\vee\cap M[p^{-1}]\right\ra$, that is there are coefficients $b_j$ such that $b_1\tilde s_1+\cdots +b_n\tilde s_n = 1$.  Arguing as above and perhaps increasing $k$, there are $\tilde b_j\in K\left\la\frac{1}{p^t}M\cap\sigma_i^\vee\right\ra$ so that:
  \[||(b_1\tilde s_1+\cdots +b_n\tilde s_n)-(\tilde b_1\tilde s_1+\cdots+\tilde b_n\tilde s_n)|| < 1.\]
  Therefore $\xi' = 1-(\tilde b_1\tilde s_1+\cdots+\tilde b_n\tilde s_n)$ is topologically nilpotent so that $1-\xi'$ is a unit in $K\left\la\sigma_i^\vee\cap p^{-t}M\right\ra$.  But it is also in the ideal generated by the $\tilde s_j$, so they generate the unit ideal.  Doing this over all the $\cU_i$, perhaps increasing $t$ a finite amount of times, we see that the $\tilde s_j$ generate $\cM^{p^t}$.
\end{proof}
Using this we can promote a well known vanishing theorem to the perfectoid setting.
\begin{Th}[Demazure Vanishing in the Perfectoid Setting]
	If $\cL$ is a globally generated line bundle on ${\X}$, then for all $i>0$
	\[\HH^i(X,\cL)=0.\]
\end{Th}
\begin{proof}
	By Theorem \ref{main2} we can find some line bundle $\cM$ on $X$ such that $\cL$ is identified with $\cM^{1/p^k}$, and Proposition \ref{globallyGenerated} shows that $\cM^{p^t}$ is globally generated for $t>>0$.  Applying Demazure vanishing \cite[Theorem 9.2.3]{CLS} we see that $\HH^i(X,\cM^{p^t}) = 0$ for all $i>0$ and $t>>0$, so that taking completed direct limits and applying Theorem \ref{cohomologyTheorem} gives the result.
\end{proof}
The proof of Batyrev-Borisov vanishing will be essentially identical, but the statement requires a bit of setup, and we will summarize without proof the necessary results.  The results are carefully described and proven over the complex numbers in \cite[Sections 3 and 4]{CLS}, and in general in \cite{Demazure} and \cite{Danilov}.

Fix a fan $\Sigma$ and let $\Sigma(1)$ be the 1 dimensional cones of $\Sigma$.  To each ray $\rho\in\Sigma(1)$ the \textit{orbit cone correspondence} \cite[Theorem 3.2.6]{CLS} canonically assigns a divisor $D_\rho\subseteq X_\Sigma$.  We also assign to $\rho$ its minimal generator in $N$, which we call $u_\rho$.  With this data, we can now assign to every $m\in M$ the divisor
\[\operatorname{div}(m) = \sum_{\rho\in\Sigma(1)}\la m,u_\rho\ra D_\rho.\]
The divisor $\operatorname{div}(m)$ is the principal divisor associated to the character $\chi^m$, and the map $\operatorname{div}$ fits into the following exact sequence \cite[Theorem 4.1.3]{CLS}
\begin{equation}\label{divisorSequence}
	M\longto\bigoplus_{\rho\in\Sigma(1)}\bZ\cdot D_\rho\longto\Cl(X)\longto0.
\end{equation}
The term in the middle is the set of Weil divisors of $X_\Sigma$ invariant under the torus action.  Central to the statement of Batyrev-Borisov vanishing is a polyhedron associated to an torus invariant divisor.
\begin{Def}[{\cite[(4.3.2)]{CLS}}]\label{polytopeDef}
	Fix a torus invariant divisor $D = \sum a_\rho D_\rho$.  Its associated polytope is
	\[P_D = \{m\in M_\R:\la m,u_\rho\ra\ge-a_\rho\text{ for all }\rho\in\Sigma(1)\}.\]
\end{Def}
\begin{Def}
	Let $S\subseteq M_{\mathbb{R}}$.  The \textit{affine hull} of $S$ is
	\[\operatorname{Aff}(S):=\left\{\sum_{i=1}^k \alpha_ix_i:x_i\in S\text{ and }\sum\alpha_i = 1\right\}.\]
	The relative interior $\operatorname{Relint}(S)$ is the interior of the affine hull of $S$.
\end{Def}
We can now state Batyrev-Borisov vanishing, first established in \cite{BB}.
\begin{Th}[Batyrev-Borisov Vanishing {\cite[Theorem 9.3.5]{CLS}}]\label{BBV}
Let $X = X_\Sigma$ be a complete toric variety over a field $k$, and $D$ a basepoint free divisor.  Then
\begin{itemize}
	\item{$\HH^i(X,\cO(-D)) = 0$ for all $i\not=\dim P_D$.}
	\item{If $i=\dim P_D$ then
	\[\HH^i(X,\cO(-D)) = \bigoplus_{m\in\operatorname{Relint}(P_D)\cap M}k\cdot\chi^{-m}.\]
	}
\end{itemize}
\end{Th}
On our perfectoid space we have not yet developed a good notion of divisors, so we need a dictionary between line bundles and divisors.  This isn't too difficult, because as $X$ is a normal noetherian scheme, there is a natural injection $\Pic X\into\Cl X$ (see for example \cite[14.2.7]{Vakil}).  Given a line bundle $\cL\in\Pic X$ we can use this injection and the exact sequence \ref{divisorSequence} to build a divisor $D = \sum a_\rho D_\rho$ such that $\cL\cong\cO(D)$, and this divisor is well defined up to an element of $M$.  Notice that essentially by definition $\cL$ is globally generated if and only if $D$ is basepoint free.  Furthermore, we can associate $\cL^n$ to the divisor $nD$.	The dimension of the polytope from Defintion \ref{polytopeDef} isolated the important cohomological degree in the Batyrev-Borisov vanishing theorem, and we can now access that integer in the perfectoid setting.
\begin{Def}\label{setup}
	Let $\Sigma$ be a complete fan, $X = X_\Sigma$ and ${\X}\to X$ the associated perfectoid cover.  Assume that the Picard group of the rigid space $\cX_\Sigma$ trivializes on the affinoid cover induced by $\Sigma$.  Fix $\cL\in\Pic{\X}$ and use Theorem \ref{main2} to identify it with $\cM^{1/p^k}$ for $\cM\in\Pic X$.  Associate to $\cM$ a divisor $D$ using sequence \ref{divisorSequence}, and consider the polytope $P_D$.  We define
	\[d_\cL:=\dim P_D.\]
\end{Def}
\begin{Lemma}\label{dLwellDef}
	The integer $d_\cL$ is well defined.
\end{Lemma}
\begin{proof}
	Suppose $D$ and $D'$ are two torus invariant divisors associated to $\cM$.  the sequence \ref{divisorSequence} asserts that $D$ and $D'$ differ by the divisor of some $m\in M$, so that by \cite[Section 4.3]{CLS}(1) $P_{D'}$ is a translate of $P_D$ by $m$, and hence share dimension.  The choice of $\cM$ is well defined up to a power of $p$, which may replace $D$ with $p^t D$.  Again by \cite[Section 4.3]{CLS} $P_{p^tD} = p^tP_D$ is a scaling and hence the dimension is unchanged.
\end{proof}
We can now prove the result.
\begin{Th}[Batyrev-Borisov Vanishing in the Perfectoid Setting]
	Consider the setup as in Definition \ref{setup}, and suppose $\cL$ is globally generated.  For all $i\not=d_\cL$:
	\[\HH^i({\X},\cL^{-1})=0\]
	If $i=d_\cL$ we have an isomorphism
	\[\HH^{d_\cL}({\X},\cL^{-1})\cong\left(\colim_n\left(\bigoplus_{m\in\operatorname{Relint}(p^n P_D)\cap M}K\cdot\chi^{-m}\right)\right)^{\wedge{}}\]
	where the transition maps are induced by $\chi^m\mapsto\chi^{pm}$.
\end{Th}
\begin{proof}
	For each $n$, we have $\cM^{-p^n} = \cO(-p^nD)$, and as in Lemma \ref{dLwellDef} the dimension of $P_{p^nD}$ is $d_\cL$.   As $\cL$ is globally generated, Proposition \ref{globallyGenerated} implies that $\cM^{p^t}$ is globally generated for $t$ large enough, and therefore for large $t$ the divisor $p^tD$ is basepoint free. The results now follow, applying Theorem \ref{BBV} to all large enough $p$-powers of $\cM$, passing to the completed direct limit, and applying Theorem \ref{cohomologyTheorem}.
\end{proof}


\begin{thebibliography}{1}

	\bibitem{BB}Batyrev V. and Borisov, L, \textit{On Calabi-Yau complete intersection in toric varieties}, Higher-dimensional complex varieties (Trento, 1994), 1996, 39-65.
	\bibitem{Peteretal}Blakestad, C., Gvirtz, D., Heuer, B., Shchedrina, D., Shimizu, K., Wear, P., and Yao, Z., \textit{Perfectoid covers of abelian varieties}, arXiv preprint arXiv:1804.04455, 2018.
	\bibitem{CaraianiScholze}Caraiani, A., and Scholze, P, \textit{On the generic part of the cohomology of compact unitary Shimura varieties}, Annals of Mathematics, 186(3) 649-766, 2017.
	\bibitem{CLS}Cox, D., Little, J., and Schenck, H., \textit{Toric varieties}, American Mathematical Soc., 2011.
	\bibitem{Danilov}Danilov, V.I., \textit{The Geometry of Toric Varieties}, Russian Math Surveys 33:2, 97-154, 1978.
  \bibitem{Demazure}Demazure, M, \textit{Sous-groupes alg\'ebriques de rang maximum du groupe de Cremona}, Annales scientifiques de l'\'E.N.S., 4e.
	\bibitem{DHthesis}Dorfsman-Hopkins, G., \textit{Projective Geometry for Perfectoid Spaces}, {M\"unster Journal of Mathematics, \emph{in press} (2021)}.
	\bibitem{GaRo} Gabber, O., and Ramero, L. Almost Ring Theory. Springer-Verlag, 2003.
	\bibitem{Gub}Gubeladze, I. D. "Anderson's Conjecture and the Maximal Monoid Class over which Projective Modules are Free." Math USSR Sb 63 165 (1989)
	\bibitem{Ben} Heuer, B., \textit{Line bundles on perfectoid covers: case of good reduction}, arXiv preprint arXiv:2105.05230, 2021.
  \bibitem{Huber1}Huber, R., \textit{A generalization of formal schemes and rigid analytic varieties}, Math. Z. 217(4) 513-555.
  \bibitem{Huber2}Huber, R., \textit{Continuous valuations}, Math. Z. 212(3) 455-477.
	\bibitem{toroidal}G. Kempf, Finn Faye Knudsen, D. Mumford, and B. Saint-Donat. \textit{Toroidal Embeddings. I.} Lecture Notes in Mathematics, Vol. 339. Springer-Verlag, Berlin, 1973.
	\bibitem{Lutke}L\"utkebohmert, W. Vektorraumb\"undel \"uber nichtarchimedischen holomorphen R\"aumen. Math. Z. 152, 2 (1977), 127–143.
  \bibitem{ScholzeDiamonds}Scholze, P. \textit{\'{E}tale Cohomology of Diamonds}, arXiv preprint arXiv:1709.07343, 2017.
	\bibitem{ScholzePadicHodgeTheory}Scholze, P., \textit{$p$-adic Hodge theory for rigid-analytic varieties}, Forum of Mathematics, Pi, 1, E1. doi:10.1017/fmp.2013.1, 2013.
	\bibitem{ScholzePS}Scholze, P, \textit{Perfectoid spaces}, Publications math\'ematiques de l IH\'ES 116.1 (2012): 245-313.
  \bibitem{Stacks}Stacks Project Authors, \textit{Stacks Project}, 2020.
	\bibitem{Vakil}Vakil, R., \textit{Foundations of Algebraic Geometry}, http://math.stanford.edu/~vakil/216blog/FOAGnov1817public.pdf, 2017.

\end{thebibliography}
\end{document}